\documentclass[a4paper]{amsart}
\usepackage{amsmath,amssymb,amsfonts}
\usepackage{amsthm}
\usepackage{mathtools}
\usepackage{paralist}
\usepackage{hyperref}
\usepackage[shortcuts]{extdash} %Correct hyphenation for words with hyphen

\hypersetup{
  colorlinks   = true, %Colours links instead of ugly boxes
  urlcolor     = blue, %Colour for external hyperlinks
  linkcolor    = blue, %Colour of internal links
  citecolor    = blue  %Colour of citations
}

\newtheorem{theorem}{Theorem}[section]
\newtheorem{algorithm}[theorem]{Algorithm}
\newtheorem{lemma}[theorem]{Lemma}
\newtheorem{corollary}[theorem]{Corollary}
\newtheorem{proposition}[theorem]{Proposition}

\theoremstyle{definition}
\newtheorem{definition}[theorem]{Definition}
\newtheorem{example}[theorem]{Example}
\newtheorem{remark}[theorem]{Remark}

\DeclareMathOperator{\Bin}{Bin}

\DeclareMathOperator{\LC}{LC}
\DeclareMathOperator{\LT}{LT}

\DeclareMathOperator{\Mat}{Mat}
\DeclareMathOperator{\lcm}{lcm}
\DeclareMathOperator{\Syz}{Syz}
\DeclareMathOperator{\Mon}{Mon}
\DeclareMathOperator{\Rad}{Rad}
\DeclareMathOperator{\mysep}{sep}
\DeclareMathOperator{\nil}{nil}
\DeclareMathOperator{\myexp}{exp}
\DeclareMathOperator{\Id}{Id}

\DeclareMathOperator{\mult}{mult}
\DeclareMathOperator{\Supp}{Supp}
\DeclareMathOperator{\vspan}{vspan}
\newcommand\sprod{\textstyle\prod\limits}
\newcommand\ssum{\textstyle\sum\limits}

\let\rho=\varrho
\let\epsilon=\varepsilon

\title{Computing the Binomial Part of a Polynomial Ideal}

\author{Martin Kreuzer}
\address[Martin Kreuzer]{Fakult\"{a}t f\"{u}r Informatik und Mathematik \\
Universit\"{a}t Passau, D-94032 Passau, Germany}
\email{martin.kreuzer@uni-passau.de}

\author{Florian Walsh}
\address[Florian Walsh]{Fakult\"{a}t f\"{u}r Informatik und Mathematik \\
Universit\"{a}t Passau, D-94032 Passau, Germany}
\email{florian.walsh@uni-passau.de}

\date{\today}

\begin{document}

\begin{abstract}
Given an ideal $I$ in a polynomial ring $K[x_1,\dots,x_n]$ over a field~$K$, we present a complete
algorithm to compute the binomial part of $I$, i.e., the subideal ${\rm Bin}(I)$ of $I$ generated by all monomials 
and binomials in $I$. This is achieved step-by-step. First we collect and extend several algorithms 
for computing exponent lattices in different kinds of fields. Then we generalize them to compute exponent
lattices of units in 0-dimensional $K$-algebras, where we have to generalize the computation of the
separable part of an algebra to non-perfect fields in characteristic $p$. Next we examine the computation
of unit lattices in affine $K$-algebras, as well as their associated characters and lattice ideals.
This allows us to calculate ${\rm Bin}(I)$ when $I$ is saturated with respect to the indeterminates
by reducing the task to the 0-dimensional case. Finally, we treat the computation of ${\rm Bin}(I)$
for general ideals by computing their cellular decomposition and dealing with finitely many special 
ideals called $(s,t)$-binomial parts. All algorithms have been implemented in {\tt SageMath}.
\end{abstract}

\keywords{binomial part, binomial ideal, exponent lattice, cellular decomposition}

\subjclass[2020]{Primary 13P05; Secondary 12-08, 13C13, 13F65.}

\maketitle

%\tableofcontents

%%%%%%%%%%%%%%%%%%%%%%%%%%%%%%%%%%%%%
%
%  Section 1: Introduction
%
%%%%%%%%%%%%%%%%%%%%%%%%%%%%%%%%%%%%%

\section{Introduction}

Let $P=K[x_1,\dots,x_n]$ be a polynomial ring over a field~$K$. Polynomials of the form
$t= x_1^{\alpha_1} \cdots x_n^{\alpha_n}$ with $\alpha_i \ge 0$ will be called {\it terms},
the set of all terms is denoted by $\mathbb{T}^n$, polynomials of the form $at$ with 
$a\in K\setminus \{0\}$ and $t\in \mathbb{T}^n$ will be called {\it monomials}, polynomials of the form 
$as + bt$ with $a,b\in K\setminus \{0\}$ and $s,t\in\mathbb{T}^n$ will be called {\it binomials},
and polynomials of the form $s-t$ with $s,t\in \mathbb{T}^n $ will be called {\it unitary binomials}.

An ideal~$I$ in~$P$ is called a {\it binomial ideal}\/ if it is generated by monomials and binomials.
These ideals are well-studied and occur in different contexts (see for
instance~\cite{eisenbud1996binomial,herzog2018binomial}). It is therefore a natural problem to search 
for binomials within a given polynomial ideal. More precisely, the ideal $\Bin(I)$ generated by 
all monomials and binomials in~$I$ is called the \textit{binomial part} of~$I$. The main topic of this paper
is to develop a general algorithm for computing generators of $\Bin(I)$.

Applications of this algorithm include a variety of problems, e.g., computing algebraic relations of 
C-finite sequences (see~\cite{kauers2008computing}), solving the constructive membership problem 
for commutative matrix groups (see~\cite{babai1996multiplicative}), and computing the Zariski closure 
of a matrix group (see~\cite{derksen2005quantum}).

To the best of our knowledge, our main result yields the first general algorithm for computing 
the binomial part of an arbitrary polynomial ideal. However, some special cases and related topics 
have been studied before. The monomial part of a polynomial ideal~$I$, i.e., the
ideal generated by all monomials contained in~$I$, is a subideal of $\Bin(I)$ and can be computed using 
homogenization (see Tutorial~50 in~\cite{kreuzer2005computational}). In~\cite{katthan2019polynomial} 
the authors construct an algorithm for checking whether an ideal is binomial after applying an ambient automorphism. 
A method for finding sparse polynomials which vanish on an algebraic set is proposed in~\cite{hauenstein2022binomiality}. 
For univariate polynomial ideals, computing the binomial part means computing a binomial multiple of lowest degree 
of its generator. Effective methods for this task are presented in~\cite{giesbrecht2012computing}. 
The computation of the binomial part of a principal ideal can be reduced to the univariate case. 

For an ideal~$I$ in $\mathbb{Q}[x_1, \dots, x_n]$ satisfying $I : \langle x_1 \cdots x_n \rangle = I$, 
an algorithm which computes $\Bin(I)$ is presented in~\cite{jensen2017finding}. It uses methods 
from tropical geometry to reduce the problem to 0-dimensional ideals. The authors also provide the following 
example which shows that a degree bound for the generators of $\Bin(I)$ would need to depend on the 
coefficients of the generators of~$I$. No such bound seems to be known. 

\begin{example}
For $n \in \mathbb{N}$, let $I \subseteq \mathbb{Q}[x,y,z]$ be the ideal generated by $(x-z)^2$ and
$nx-y-(n-1)z$. Then $x^n-yz^{n-1}$ is the binomial of least degree contained in~$I$.
\end{example}

To achieve our main goal of constructing an algorithm for computing the binomial part
of an arbitrary ideal~$I$ in~$P$, we use a chain of reductions. These are
presented from the bottom up, starting with the case where~$I$ is a maximal ideal.

The initial step in Section~2 is the computation of exponent lattices in fields. Given elements
$f_1, \dots, f_k$ in a multiplicative abelian group, for instance the group of units of a ring,
the set of all $(a_1,\dots,a_k) \in \mathbb{Z}^k$ such that $f_1^{a_1}\cdots f_k^{a_k} = 1$
forms a lattice in~$\mathbb{Z}^k$ which is called the {\it exponent lattice} of $(f_1,\dots,f_k)$.
Algorithms for computing exponent lattices in various fields are known (see for example~\cite{derksen2005quantum}, 
\cite{ge1993algorithms}, \cite{kauers2005algorithms}, and~\cite{zheng2019effective}). Here we provide
variants of these algorithms in a consistent notation which focus on a compact presentation and ease of implementation.
Why are we interested in such exponent lattices? For a maximal ideal~$I$ in a polynomial ring $P=K[x_1,\dots,x_n]$
over a field~$K$ such that $x_i \notin I$ for $i=1,\dots,n$, the unitary binomials in~$I$ correspond 1--1
to the elements of the exponent lattice of $(\bar{x}_1,\dots,\bar{x}_n)$ in the field $P/I$.

Next, in Section~3, we generalize these algorithms and show how exponent lattices of units
in 0-dimensional affine $K$-algebras can be computed.
Over the base field $K=\mathbb{Q}$, exponent lattices can be computed using the algorithm presented
in Section~8 of~\cite{lenstra2018algorithms} for which we present a slightly generalized version
(see~\ref{alg:zero_dim_lattice0}). The main result in this section is
a new algorithm for 0-dimensional algebras in finite characteristic (see Algorithm~\ref{alg:zero-dimP}).
It is based on the fact that the algebra can be split effectively into a direct
sum of its separable part and its nilradical (see Algorithm~\ref{alg:sep_dec} which generalizes Algorithm~5.5.6
in~\cite{kreuzer2016computational}). More precisely, we show that such a decomposition can also
be achieved over non-perfect fields. 

In the remaining sections we generalize these algorithms even further as follows. Given a ring~$R$, a tuple of elements
$F=(f_1,\dots,f_k) \in R^k$, and a subgroup~$G$ of the group of units $R^\times$, the lattice consisting
of all $a=(a_1,\dots,a_k) \in \mathbb{Z}^k$ such that
$$
f_1^{a^+_1} \cdots f_k^{a^+_k} \;-\; g\cdot f_1^{a^-_1} \cdots f_k^{a^-_k} = 0 \quad\hbox{\rm for some\ }g\in G
$$
is called the \textit{unit lattice} of $F=(f_1,\dots,f_k)$ with respect to~$G$. 
Here we let $a=a^+ - a^- = (a^+_1 - a^-_1,\dots, a^+_k-a^-_k)$ be the unique decomposition
with $a^+_i = \max\{a_i,0\}$ and $a^-_i = \min\{a_i,0\}$.
Thus the exponent lattice of a tuple~$F$ is nothing but its unit lattice with respect to $G=\{1\}$.
For us, the most important rings for which we want to compute unit lattices are affine $K$-algebras
$R=P/I$. Namely, if we use the tuple of residue classes $(\bar{x}_1,\dots,\bar{x}_n) \in R^n$ and the group 
$G=K^\times$, the elements of the unit lattice correspond to the binomials in~$I$.
 
Thus we start to examine the computation of unit lattices of affine algebras $R=P/I$ in Section~4.
The first case we consider is the case when~$I$ is saturated with respect to the tuple $F=(f_1,\dots,f_k)$,
i.e., when $I :_P \langle f_1\cdots f_k\rangle = I$. In this case the tuple give rise to a 
well-defined group homomorphism $\rho:\; \Lambda \longrightarrow G$ from its unit lattice~$\Lambda$
to the group~$G$ which is called its \textit{associated character} (see Proposition~\ref{prop:unit_lattice}).
Conversely, a lattice $\Lambda$ in~$\mathbb{Z}^n$ together with a character $\rho:\; \Lambda \longrightarrow
K^\times$ yields a binomial ideal 
$$
I_{\Lambda,\rho} \;=\; \langle x_1^{a^+_1} \cdots x_n^{a^+_n} \;-\; \rho(a) \cdot x_1^{a^-_1} \cdots
x_n^{a^-_n} \;\mid\; a = a^+ - a^- \in \Lambda \rangle 
$$ 
in~$P$ which is called the \textit{lattice ideal} associated to $(\Lambda,\rho)$.
In particular, we check that if a given ideal~$I$ in~$P$ satisfies $I :_P \langle x_1\cdots x_n \rangle =I$, 
then $\Bin(I)$ is a lattice ideal (see Corollary~\ref{cor:lattice_ideal}).
In Algorithm~\ref{alg:lattice_intersection} we discuss a method for computing intersections of lattices with
associated characters, in Proposition~\ref{prop:invertible} we note that lattices with associated characters localize,
and in Proposition~\ref{prop:field_extension} we verify that they are stable under base field extensions.

In Section~5 we start the actual computation of unit lattices of tuples in an affine $K$-algebra
$R=P/I$ with respect to the group~$K^\times$. By Corollary~\ref{cor:lattice_ideal}, this yields
an algorithm for computing $\Bin(I)$ when~$I$ is saturated with respect to $x_1\cdots x_n$.
For ideals in $\mathbb{Q}[x_1, \dots, x_n]$, such an algorithm was formulated in~\cite{jensen2017finding}. 
We present an alternative approach. It avoids the use of methods from tropical geometry, for which there seems to exist no
implementation so far. Based on the computation of a maximal set of independent indeterminates, we first reduce the task to
the case of 0-dimensional ideals~$I$ (see Algorithm~\ref{alg:unit_lattice}).
Then we reduce the task of computing unit lattices in~$R$ to exponent lattices in~$R$ 
(see Proposition~\ref{prop:reduction_to_unitary} and Algorithm~\ref{alg:partial_char_zero_dim}).

Finally, in Section~6, we reduce the computation of the binomial part of an ideal~$I$ in~$P$
to the case $I :_P \langle x_1\cdots x_n\rangle =I$. The problem is that if an ideal~$I$ in~$P$ 
does not satisfy this condition, then the binomials in~$I$ are no longer in correspondence with a unit lattice.
For the desired reduction, we employ the decomposition of~$I$ into cellular ideals 
(see Algorithm~\ref{alg:cellular_decomposition}). Here an ideal ~$I$ is called \textit{$Y$-cellular}
for a set~$Y$ of indeterminates in $\{x_1,\dots,x_n\}$ if~$I$ is saturated with respect to the 
indeterminates in~$Y$ and the remaining indeterminates are nilpotent modulo~$I$.
For a $Y$-cellular ideal, the main task in computing $\Bin(I)$ can be reduced to computing
$(s,t)$-binomial parts, where $s,t$ are terms in $K[X\setminus Y]$. This is achieved in
Algorithm~\ref{alg:binomials_fixed} and allows us to compute the binomial part of a cellular ideal.
The final step is to find the binomial part of an intersection of cellular ideals
which is done in Algorithm~\ref{alg:binomial_part}. The paper concludes with two optimizations:
a restriction of the set of pairs $(s,t)$ which has to be considered (using Algorithm~\ref{alg:vr_binomials}), 
and a simplification of the entire algorithm in the case of a radical ideal~$I$ (see Algorithm~\ref{alg:bin_of_rad}).

The definitions and notation in this paper follow the books~\cite{kreuzer2000computational} and~\cite{kreuzer2005computational}.
An important aspect is that all algorithms in this paper have been implemented using the software system
\texttt{SageMath}~\cite{sagemath}. The complete package is available freely from the second author's GitHub
page~\footnote{{\tt https://github.com/abacus42/binomial-part}}

\bigskip\bigbreak
%%%%%%%%%%%%%%%%%%%%%%%%%%%%%%%%%%%%%%%%%%%%%
%
%  Section 2: Exponent Lattices in Fields
%
%%%%%%%%%%%%%%%%%%%%%%%%%%%%%%%%%%%%%%%%%%%%%

\section{Exponent Lattices in Fields}\label{sec2}

Let $K$ be a field, and let $\mathfrak{m}$ be a maximal ideal in $K[x_1,\dots,x_n]$ such that 
$x_i \notin \mathfrak{m}$ for $i=1,\dots,n$. Then $L=P/\mathfrak{m}$ is a field, 
and the residue classes $\overline{x}_i$ of $x_i$ in $L$ are units. 
In this setting, the unitary binomials in~$\mathfrak{m}$ are determined by all 
$a =(a_1,\dots,a_n) \in \mathbb{Z}^n$ with $\overline{x}_1^{a_1} \cdots \overline{x}_n^{a_n} = 1$
in~$L$. This motivates the following definition.

\begin{definition}
Let $G$ be a multiplicative abelian group and $f_1,\dots,f_k \in G$. Then the lattice
$$
\Lambda = \{(a_1,\dots,a_k) \in \mathbb{Z}^k \mid f_1^{a_1} \cdots f_k^{a_k} = 1 \}
$$ 
is called the \textbf{exponent lattice} of $(f_1,\dots,f_k)$. 
    
If $G$ is the group of units of a ring~$R$, we also refer to~$\Lambda$ as the \textbf{exponent lattice} 
of $(f_1,\dots,f_k)$ in~$R$.
\end{definition}

This section provides an overview on how exponent lattices in different types of fields can be computed. Algorithms for
this problem are already known but scattered throughout the literature. We present variants of these algorithms,
which allow a straightforward implementation in a computer algebra system such as \texttt{SageMath}.
In subsequent sections we show how the computation of the binomial part of a polynomial ideal can be
reduced to computing exponent lattices in fields. Let us begin with the case of finite fields.

\begin{algorithm}{\bf (Computing Exponent Lattices in Finite Fields)}\label{alg:finite_field}\\
Let $q$ be a prime power, let $K = \mathbb{F}_q$ be the finite field with~$q$ elements, and let $f_1,\dots,f_k \in K^\times$.
The following sequence of instructions forms an algorithm which computes the exponent lattice of $(f_1,\dots,f_k)$.
\begin{enumerate}
\item[(1)] Compute a generator $g$ of the cyclic group $K^\times$.

\item[(2)] For $i=1,\dots,k$, compute the discrete logarithms $e_i:=\log_g(f_i)$.

\item[(3)] Compute the solution space $\Lambda' \subseteq \mathbb{Z}^{k+1}$ of the linear equation over the integers in
the indeterminates $y_1, \dots, y_{k+1}$ given by
$$
e_1 y_1 + \cdots + e_k y_k + y_{k+1}(q-1) = 0.
$$

\item[(4)] Return the projection~$\Lambda$ of the solution space $\Lambda'$ onto its first~$k$ coordinates.

\end{enumerate}
\end{algorithm}

\begin{proof}
A tuple $(a_1,\dots,a_k) \in \mathbb{Z}^k$ is in the exponent lattice of $(f_1, \dots, f_k)$ if and only if
$g^{e_1^{a_1}} \cdots g^{e_k^{a_k}} = 1$ in~$K$. This is the case if and only if the order $q-1$ of~$K^\times$ divides
$e_1^{a_1}+ \cdots+ e_k^{a_k}$, which is equivalent to $(a_1,\dots,a_k)$ being a projection of an element of $\Lambda'$ 
onto its first~$k$ coordinates.
\end{proof}

For number fields, Masser~\cite{masser1988linear} published a bound for the norm of the basis elements
of an exponent lattice. This yields an algorithm based on exhaustive search. Different and more efficient methods are
presented by Kauers in~\cite{kauers2005algorithms} and by Zheng and Xia in~\cite{zheng2019effective}. A polynomial time
algorithm was developed by Ge in his PhD thesis~\cite{ge1993algorithms}.
Let us present an approach which works by first determining the integral unit lattice of the given elements which we now define. 
The ring of integers in a number field~$K$ is denoted by~$\mathcal{O}_K$.

\begin{definition}
Let $K$ be a number field, and let $f_1, \dots, f_k \in K^\times$.
The lattice $\{(a_1,\dots,a_k) \in \mathbb{Z}^k \mid f_1^{a_1} \cdots f_k^{a_k} \in \mathcal{O}_K^\times\}$ is called the
\textbf{integral unit lattice} of $(f_1, \dots, f_k)$.
\end{definition}

Integral unit lattices are special kinds of unit lattices, as defined later.
To compute integral unit lattices, we use the fact that every fractional ideal in~$K$ can be written as a product of
non-zero prime ideals in~$\mathcal{O}_K$ and their inverses. For a prime ideal~$\mathfrak{p}$, we denote the multiplicity 
with which $\mathfrak{p}$ occurs in the prime factorization of a fractional ideal~$I$ by $\mult_\mathfrak{p}(I)\in \mathbb{Z}$.

\begin{algorithm}{\bf (Computing Integral Unit Lattices)}\label{alg:integral_unit_lattice}\\
Let $K$ be a number field. The following sequence of instructions forms an algorithm which computes the integral unit
lattice of $(f_1,\dots,f_k) \in (K^\times)^k$.
\begin{enumerate}
\item[(1)] For $i=1,\dots,k$, form the fractional ideal generated by~$f_i$ and compute its factorization into
prime ideals. Let $\{\mathfrak{p}_1, \dots, \mathfrak{p}_m\}$ be the set of all prime ideals occurring
in these factorizations.

\item[(2)] Return the solution space $\Lambda \subseteq \mathbb{Z}^k$ of the linear system of equations over~$\mathbb{Z}$
in the indeterminates $y_1, \dots, y_k$ given by
$$
\ssum_{i=1}^k \mult_{\mathfrak{p}_j}(\langle f_i \rangle) y_i = 0\qquad \hbox{\it for\ } j=1,\dots,m.
$$

\end{enumerate}
\end{algorithm}

\begin{proof}
The fractional ideals in~$K$ form a free abelian group generated by the non-zero prime ideals in~$\mathcal{O}_K$.
Consequently, a tuple $(a_1,\dots,a_k) \in \mathbb{Z}^k$ is in the integral unit lattice of $(f_1,\dots,f_k)$ if and
only if
$$
\mult_{p_j}(f_1^{a_1} \cdots f_k^{a_k}) = a_1 \mult_{p_j}(f_1) + \cdots + a_k \mult_{p_j}(f_k) = 0
$$
for all $j=1,\dots,m$. This is equivalent to $(a_1,\dots,a_k) \in \Lambda$.
\end{proof}

For $K = \mathbb{Q}$, we can write each element $f_i$ of the input of this algorithm as a fraction $g_i/h_i$ with
$g_i, h_i$ in~$\mathbb{Z} \setminus \{0 \}$. Step~(1) then simplifies to determining a set of pairwise coprime
integers such that each of the integers $g_i$ and $h_i$ can be written as a product of elements in this set. Such a set
of pairwise coprime integers can be computed in essentially linear time (see~\cite{bernstein2005factoring}). Also for
general number fields, the factorizations into prime ideals in Step~(1) of this algorithm can be avoided. In fact,
Ge~\cite{ge1994recognizing} provides an algorithm which computes the integral unit lattice in polynomial time.

After computing a basis $b_1, \dots, b_m$ with $b_i = (b_{i1}, \dots, b_{ik}) \in \mathbb{Z}^k$ of the integral unit
lattice of a tuple $(f_1,\dots,f_k)$ in~$(K^\times)^k$, we form the elements $g_i = f_1^{b_{i1}}\cdots f_k^{b_{ik}}$
in~$\mathcal{O}_K^\times$. The exponent lattice of $(f_1,\dots,f_k)$ can then be determined by computing the exponent
lattice of $(g_1,\dots,g_m)$ in~$(\mathcal{O}_K^\times)^m$. By Dirichlet's unit theorem, we know that
$\mathcal{O}_K^\times$ is a finitely generated abelian group with a cyclic torsion subgroup. This yields the
following algorithm.

\begin{algorithm}{\bf (Computing Exponent Lattices in Number Fields)}\label{alg:number_field}\\
Let $K$ be a number field and let $\mathcal{O}_K$ be its ring of integers. The following sequence of instructions forms
an algorithm which computes the exponent lattice of $(f_1,\dots,f_k) \in (K^\times)^k$.
\begin{enumerate}
\item[(1)] Using Algorithm~\ref{alg:integral_unit_lattice}, compute a basis $b_1, \dots, b_m \in \mathbb{Z}^k$ of the 
integral unit lattice of $(f_1,\dots,f_k)$.

\item[(2)] For $i=1,\dots,m$, form the elements $g_i = f_1^{b_{i1}}\cdots f_k^{b_{ik}} \in \mathcal{O}_K^\times$.

\item[(3)] Compute a system of fundamental units $\epsilon_1, \dots, \epsilon_\ell$ of $\mathcal{O}_K^\times$, as well as
a generator $\zeta$ of the cyclic group of roots of unity in~$K$. Let~$r$ be the order of~$\zeta$.

\item[(4)] For $i=1,\dots,m$, write $g_i = \epsilon_1^{c_{i1}} \cdots \epsilon_\ell^{c_{i\ell}} \cdot \zeta^{d_i}$
    with $c_{ij}, d_i \in \mathbb{Z}$.

\item[(5)] Compute the solution space $\Lambda' \subseteq \mathbb{Z}^{m+1}$ of the linear system of equations
over~$\mathbb{Z}$ in the indeterminates $y_1, \dots, y_{m+1}$ given by
$$
\ssum_{i=1}^m c_{ij}y_i = 0 \quad \text{\ \it for\ } j=1, \dots, \ell \quad \text{\ \it and\ } \quad
\ssum_{i=1}^m d_i y_i = r y_{i+1}.
$$

\item[(6)] Let $\Lambda$ be the projection of $\Lambda'$ onto the first~$m$ components. Return the lattice
$$
\{ (c_1b_{11}+\cdots + c_m b_{m 1}, \dots, c_1b_{1 \ell}+\cdots + c_m b_{m \ell}) \mid (c_1,\dots,c_m) \in \Lambda \}
\subseteq \mathbb{Z}^k.
$$
\end{enumerate}
\end{algorithm}

\begin{proof}
By Dirichlet's unit theorem, we know that the group $\mathcal{O}_K^\times$ is the direct product of the 
free abelian group generated by the fundamental units and the group of roots of unity in~$K$. 
Hence, $(a_1, \dots, a_m) \in \mathbb{Z}^m$ is an element of the exponent lattice of $(g_1,\dots,g_m)$ 
in $(\mathcal{O}_K^\times)^m$ if and only if
$$
\epsilon_1^{a_1c_{11}+ \cdots + a_m c_{m1}} \cdots \epsilon_\ell^{a_1c_{1\ell}+ \cdots + a_m c_{m\ell}} = 1 \quad
\text{and} \quad \zeta^{a_1d_1 + \cdots + a_m d_m} = 1.
$$
This is equivalent to $(a_1,\dots,a_m) \in \Lambda$. Step~(6) therefore returns the exponent lattice of $(f_1,\dots,f_k)$.
\end{proof}

In Algorithm~\ref{alg:partial_char_zero_dim} it will be necessary to compute exponent lattices in fields of the form
$K(x_1, \dots, x_n)$, where $x_1, \dots, x_n$ are indeterminates. If exponent lattices in the base field $K$ can be
effectively computed, then this can be achieved as follows. For a monic irreducible
polynomial~$p$ we denote the multiplicity with which~$p$ occurs in the factorization of a polynomial~$f$ by~$\mult_p(f)$.

\begin{algorithm}{\textbf{(Computing Exponent Lattices in Function Fields)}}\\\label{alg:exp_lattice_function_field}
Let $K$ be a field in which exponent lattices can be effectively computed, let $x_1,\dots,x_n$ be indeterminates, and
let $f_1, \dots, f_k \in K(x_1, \dots, x_n)$. The following sequence of instructions forms an algorithm which computes
the exponent lattice of $(f_1,\dots,f_k)$.
\begin{enumerate}
\item[(1)] For $i=1, \dots, k$ write $f_i = \lambda_i g_i/h_i$ with $\lambda_i \in K^\times$ and monic polynomials
            $g_i, h_i$ in~$K[x_1, \dots, x_n] \setminus \{0\}$.

\item[(2)] Compute the exponent lattice $\Lambda \subseteq \mathbb{Z}^k$ of $(\lambda_1, \dots, \lambda_k)$ in $K^\times$.

\item[(3)] For $i=1, \dots, k$ compute the factorizations of $g_i$ and $h_i$ into monic irreducible polynomials. Let
            $P = \{p_1, \dots, p_m\}$ be the set of all irreducible polynomials occurring in these factorizations.

\item[(4)] Compute the solution space $M \subseteq \mathbb{Z}^k$ of the linear system of equations over~$\mathbb{Z}$
            in the indeterminates $y_1, \dots, y_k$ given by
            $$\ssum_{i=1}^k (\mult_{p_j}(f_i)-\mult_{p_j}(g_i))y_i = 0$$
            for $j=1, \dots, m$.

\item[(5)] Return the lattice $\Lambda \cap M$.

\end{enumerate}
\end{algorithm}

\begin{proof}
A tuple $(a_1,\dots,a_k) \in \mathbb{Z}^k$ is in the exponent lattice of $(f_1,\dots,f_k)$ if and only if
$\lambda_1^{a_1} \cdots \lambda_k^{a_k} = 1$ and $g_1^{a_1} \cdots g_k^{a_k}h_1^{-a_1} \cdots h_k^{-a_k} = 1$.
This is the case if and only if $(a_1,\dots,a_k) \in \Lambda$ and
$$
\mult_{p_j}(g_1^{a_1} \cdots g_k^{a_k}h_1^{-a_1} \cdots h_k^{-a_k}) =
    \ssum_{i=1}^k (\mult_{p_j}(g_i)-\mult_{p_j}(h_i))a_i = 0
$$
for $j=1,\dots,m$, which is equivalent to $(a_1,\dots,a_k) \in \Lambda \cap M$.
\end{proof}

Again, instead of computing the factorizations into irreducible polynomials in Step~(3) it is enough to compute
a set of pairwise coprime polynomials such that each of the polynomials $g_i$, $h_i$ can be written as a product of
elements from this set.

Let $K$ be any of the fields considered above, and let~$L$ be a finitely generated extension of~$K$. Then in general
the previous algorithms can not be applied directly. Instead we can use the algorithm sketched in Section~3.2
in~\cite{derksen2005quantum} to reduce the task to the case of a finite extension of~$K$. It uses the following
well-known facts.

\begin{lemma}\label{lemma:local_valuation}
Let $R$ be a local domain with non-zero and principal maximal ideal $\mathfrak{m} = \langle t \rangle$ 
such that $\bigcap_{n\geq 0} \mathfrak{m}^n = 0$. Let~$K$ be the fraction field of~$R$.  
Then every element $a \in K^\times$ can be written as $a = u t^n$ with $n \in \mathbb{Z}$ and with $u \in R^\times$. 
The map $\nu : K^\times \rightarrow \mathbb{Z}$ given by $\nu(a) = n$ is a discrete valuation and~$R$ 
is the discrete valuation ring of $\nu$.
\end{lemma}

\begin{proof}
See \cite{altman2013term}, Lemma 23.3.
\end{proof}

The localization of a normal domain~$R$ at a prime ideal~$\mathfrak{p}$ of height one satisfies the requirements of
Lemma~\ref{lemma:local_valuation} and is therefore a discrete valuation ring. We denote its discrete valuation by
$\nu_\mathfrak{p}$. The following lemma allows us to compute $\nu_\mathfrak{p}(f)$ for all $f \in R$.

\begin{lemma}\label{lemma:computing_valuation}
    Let $\mathfrak{p}$ be a height one prime ideal of a normal domain $R$, and let $r \geq 1$. For an element
    $f \in R$ we have $\nu_{\mathfrak{p}}(f) \geq r$ if and only if
    $(\mathfrak{p}^r:\langle f \rangle) \not\subseteq \mathfrak{p}$.
\end{lemma}

\begin{proof}
    Let the $g$ be the generator of the maximal ideal $\mathfrak{p}R_{\mathfrak{p}}$. Then we have
    $\nu_{\mathfrak{p}}(f) \geq r$ if and only if $g^r \mid f$ which is equivalent to
    $fR_{\mathfrak{p}} \subseteq \mathfrak{p}^rR_{\mathfrak{p}}$. The containment holds
    if and only if there exists $h \in R \setminus \mathfrak{p}$ such that $h \langle f \rangle \subseteq \mathfrak{p}^r$.
    This is equivalent to $(\mathfrak{p}^r : \langle f \rangle) \not\subseteq \mathfrak{p}$.
\end{proof}

Let $f_1, \dots, f_k \in L^\times$, and let~$t$ be an indeterminate. In the following
we form the integral closure $\overline{S}$ of the ring $S=K[f_1t, \dots, f_kt, t]$.

\begin{proposition}\label{prop:integral_closure}
    Let $S$ be defined as above, and let $K'$ be the integral closure of~$K$ within~$L$. Then the following holds.
    \begin{enumerate}
        \item[(a)] $\overline{S} \cap L = K'$
        \item[(b)] $f_1^{a_1} \cdots f_k^{a_k} \in (K')^\times$ if and only if for each height one prime ideal
            $\mathfrak{p}$ of $\overline{S}$ we have $\nu_{\mathfrak{p}}(f_1^{a_1} \cdots f_k^{a_k}) = 0$.
    \end{enumerate}
\end{proposition}

\begin{proof}
    Part~(a) is shown in Theorem~6.7.3 in~\cite{brennan1993effective}. For~(b), recall that $\overline{S}$ is the
    intersection of all its localizations at prime ideals of height one. We therefore have $f_1^{a_1} \cdots f_k^{a_k}$
    in $\overline{S}^\times$ if and only if its is contained in all $\overline{S}_{\mathfrak{p}}^\times$. By the
    definition of the discrete valuation in Lemma~\ref{lemma:local_valuation} this equivalent to
    $\nu_{\mathfrak{p}}(f_1^{a_1} \cdots f_k^{a_k}) = 0$ for each $\mathfrak{p}$. The claim now follows from~(a).
\end{proof}

Note that $K'$ is a finite extension of~$K$. This proposition therefore allows us to compute the exponent lattice
of $(f_1, \dots, f_k)$ in $L$ as follows.

\begin{algorithm}{\textbf{(Exponent Lattices in Finitely Generated Extensions)}}\\
Let $K$ be a perfect field such that the exponent lattice in a finite extension of $K$ can be effectively computed.
Let~$L$ be a finitely generated field extension of $K$. The following instructions form an algorithm which
computes the exponent lattice of $(f_1, \dots, f_k) \in L^k$.
\begin{enumerate}
\item[(1)] Let $t$ be an indeterminate. Compute the integral closure $\overline{S}$ of the ring
            $S=K[f_1t, \dots, f_kt, t]$.

\item[(2)] Compute the associated primes $\mathfrak{p}_1, \dots, \mathfrak{p}_s$ of the ideals $\langle t \rangle$ and
            $\langle f_it \rangle$ in $\overline{S}$ for $i=1, \dots, k$.

\item[(3)] For each $\mathfrak{p}_j$ and for each $g \in \{f_1t, \dots, f_kt, t \}$ compute the
            smallest number $r > 0$ such that $(\mathfrak{p}_j^r:\langle g \rangle) \subseteq \mathfrak{p}_j$ and
            obtain $\nu_{\mathfrak{p}_j}(g) = r-1$.

\item[(4)] For $i=1, \dots, k$ and $j=1, \dots, s$ compute
            $\nu_{\mathfrak{p}_j}(f_i) = \nu_{\mathfrak{p}_j}(f_it) - \nu_{\mathfrak{p}_j}(t)$.

\item[(5)] Compute a $\mathbb{Z}$-basis $b_1, \dots, b_\ell \in \mathbb{Z}^k$ of the solution space of the linear
            system of equations in the indeterminates $y_1, \dots, y_k$ given by
            $$y_1 \nu_{\mathfrak{p}_j}(f_1)+ \cdots + y_k \nu_{\mathfrak{p}_j}(f_k) = 0$$
            for $j=1, \dots, s$.

\item[(6)] For $i=1, \dots, \ell$ form the elements $g_i = f_1^{b_{i1}} \cdots f_k^{b_{ik}}$. Let
            $M = K[g_1, \dots, g_\ell]$ be the finite extension field of $K$ obtained by adjoining $g_1, \dots, g_\ell$.

\item[(7)] Using Algorithm~\ref{alg:number_field} or Algorithm~\ref{alg:finite_field} compute the exponent
            lattice~$\Lambda$ of $(g_1, \dots, g_\ell)$ in the field $M$.

\item[(8)] Return the lattice
$$
\{ (c_1b_{11}+\cdots + c_\ell b_{\ell 1}, \dots, c_1b_{1 k}+\cdots + c_\ell b_{\ell k}) \mid c \in \Lambda \} \subseteq \mathbb{Z}^k.
$$

\end{enumerate}
\end{algorithm}

\begin{proof}
    Firstly, we note that Steps~(2) to~(4) correctly compute the discrete valuations~$\nu_{\mathfrak{p}}(f_i)$ by
    Lemma~\ref{lemma:computing_valuation}. Since $\nu_{\mathfrak{p}}(f) = 0$ for all $\mathfrak{p}$ with
    $f \notin \mathfrak{p}$, it is enough to consider the associated primes of $\langle t \rangle$ and
    $\langle f_it \rangle$ for $i=1, \dots, k$ in Step~(2).

    If $a = (a_1, \dots, a_k) \in \mathbb{Z}^k$ is in the exponent lattice of $(f_1, \dots, f_k)$, then we have
    $f_1^{a_1} \cdots f_k^{a_k} \in (K')^\times$ where $K'$ is the integral closure of $K$ within $L$. By
    Proposition~\ref{prop:integral_closure} we have $f_1^{a_1} \cdots f_k^{a_k} \in (K')^\times$ if and only if $a$ is
    a solution of the linear system given in Step~(5). Therefore there exist $c_1, \dots, c_\ell$ such that
    $a = c_1 b_1 + \cdots + c_\ell b_\ell$. Now $a$ is in the exponent lattice of $(f_1, \dots, f_k)$ if and only if
    $(c_1, \dots, c_\ell)$ is in the exponent lattice of $(g_1, \dots, g_\ell)$. Clearly, $M$ is a subfield of $K'$.
    Therefore this exponent lattice can be effectively computed. Finally, it follows by the definition of
    $g_1, \dots, g_\ell$ that Step~(8) returns the exponent lattice of $(f_1, \dots, f_k)$.
\end{proof}

The integral closure over the perfect field $K$ in Step~(1) of this algorithm can be computed using for example the
method given in~\cite{greuel2010normalization}.

\bigskip\bigbreak
%%%%%%%%%%%%%%%%%%%%%%%%%%%%%%%%%%%%%%%%%%%%%%%%%%%%%%%%%%%%%%%%%%%%%%
%
%  Section 3: Exponent Lattices in Zero-Dimensional Affine Algebras
%
%%%%%%%%%%%%%%%%%%%%%%%%%%%%%%%%%%%%%%%%%%%%%%%%%%%%%%%%%%%%%%%%%%%%%%

\section{Exponent Lattices in Zero-Dimensional Affine Algebras}%
\label{sec3}

Let $K$ be a field. The goal of this section is to provide algorithms for computing the exponent lattice of tuples of
elements in 0-dimensional $K$-algebras. Let us begin by studying the structure of 0-dimensional affine $K$-algebras.

%%%%%%%%%%%%%%%%%%%%%%%%%%%%%%%
% Subsection 3.1
%%%%%%%%%%%%%%%%%%%%%%%%%%%%%%%
\bigbreak
\subsection*{The Structure of Zero-Dimensional Algebras}

In the following we let $P = K[x_1,\dots,x_n]$ be a polynomial ring over a field~$K$, 
and let~$I$ a 0-dimensional ideal in~$P$. Then $R = P/I$ is a 0-dimensional affine $K$-algebra. 
In particular, it is a finite dimensional $K$-vector space, and it has finitely many maximal ideals.
The intersection of all its maximal ideals is called the \textbf{zero radical} of~$R$. We denote it by $\Rad(0)$.
In Section 5.5 of~\cite{kreuzer2016computational} it is shown that if $K$ is a perfect field, then~$R$ can always
be decomposed into the direct sum of its separable subalgebra and its zero radical. Over non-perfect fields such a
decomposition need not exist. However, in the following we use ideas from~\cite{steel2005conquering} 
and extend the field until we obtain the desired decomposition.

\begin{definition}
Let $K$ be a field, let $\overline{K}$ be its algebraic closure, and let~$x$ be an indeterminate.
\begin{enumerate}
\item[(a)] A polynomial $f \in K[x]$ is called \textbf{separable} if it is either a non-zero constant 
or it factors in $\overline{K}[x]$ into pairwise distinct linear factors.

\item[(b)] Let $f=(x-a_1)^{e_1} \cdots (x-a_k)^{e_k} \in \overline{K}[x]$ with pairwise distinct 
elements $a_1, \dots, a_k\in \overline{K}$ and with $e_1, \dots, e_k \in \mathbb{N}$. 
Then $\mysep(f) = (x-a_1) \cdots (x-a_k)$ is called the \textbf{separable part} of~$f$.

\item[(c)] An element $a \in R$ of a 0-dimensional affine $K$-algebra is called \textbf{separable} 
if its minimal polynomial is separable.

\end{enumerate}
\end{definition}

Separable polynomials can be characterized as follows.

\begin{proposition}
For $f \in K[x]$, the following are equivalent.
\begin{enumerate}
\item[(a)] The polynomial $f$ is separable.
       
\item[(b)] We have $\gcd(f,f') = 1$.

\item[(c)] For every extension field~$L$ of~$K$, the polynomial~$f$ 
is squarefree in $L[x]$.

\end{enumerate}
\end{proposition}

\begin{proof}
See \cite{becker1993groebner}, Proposition 7.33.
\end{proof}

To compute a purely inseparable extension $L$ of $K$ and the separable part
$\mysep(f) \in L[x]$ of a polynomial $f \in K[x]$, we can use the algorithm 
given in Section~3 of~\cite{steel2005conquering} or Algorithm~1 in~\cite{kemper2002calculation}.

    \begin{lemma}\label{lemma:min_poly}
        Let $I$ be a 0-dimensional ideal in $K[x_1, \dots, x_n]$, let $L$ be an extension field of $K$, and let
        $a\in K[x_1, \dots, x_n]/I$. Consider the canonical $K$-algebra homomorphism
        $$\varphi: K[x_1, \dots, x_n]/I \rightarrow L[x_1, \dots, x_n]/IL[x_1, \dots, x_n].$$ The minimal polynomial
        $\mu_a$ of $a$ and the minimal polynomial $\mu_{\varphi(a)}$ of~$\varphi(a)$ coincide.
    \end{lemma}
    \begin{proof}
        The map $\varphi$ is injective since $IL[x_1, \dots, x_n] \cap K[x_1, \dots, x_n] = I$. Therefore
        $0 = \mu_{\varphi(a)}(\varphi(a)) = \varphi(\mu_{\varphi(a)}(a))$ implies $\mu_{\varphi(a)}(a) = 0$. This shows
        that $\mu_a$ divides $\mu_{\varphi(a)}$. The observation $0 = \varphi(\mu_a(a)) = \mu_a(\varphi(a))$ then
        shows $\mu_{\varphi(a)} \mid \mu_a$.
    \end{proof}

    The following is a generalization of Proposition~5.5.2 in~\cite{kreuzer2016computational}.

    \begin{proposition}\label{prop:sep_algebras}
        Let $I$ be a 0-dimensional ideal in $K[x_1, \dots, x_n]$, and let $R = K[x_1, \dots, x_n]/I$. Let
        $\overline{x}_i$ denote the residue class of~$x_i$ in~$R$. Then the following conditions are equivalent.
        \begin{enumerate}[(a)]
            \item For every extension field $L$ of $K$ the ring $R \otimes_K L$ is reduced, i.e.,
                $IL[x_1, \dots, x_n]$ is radical.
            \item The elements $\overline{x}_1, \dots, \overline{x}_n$ are separable.
            \item All elements in $R$ are separable.
        \end{enumerate}
    \end{proposition}
    \begin{proof}
        We first prove $(a) \Rightarrow (c)$. Let $a \in K[x_1, \dots, x_n]$, let $f \in K[z]$ be the minimal
        polynomial of $a+I$ in $R$, and let $L$ be an extension field of $K$. The minimal polynomial of
        $a +IL[x_1, \dots, x_n]$ in $L[x_1, \dots, x_n]/IL[x_1, \dots, x_n]$ equals~$f$, by Lemma~\ref{lemma:min_poly}.
        Let $g \in L[z]$ be the squarefree part of $f$ considered as polynomial in $L[z]$. Then
        $f(a) \in IL[x_1, \dots, x_n]$ implies that $g(a)$ is an element of
        $\Rad(IL[x_1, \dots, x_n]) = IL[x_1, \dots, x_n]$. Consequently, we have $g = f$. This shows that $f$ is
        squarefree in $L[z]$ for every extension field $L$ of $K$.

        Clearly, (c) implies (b). The remaining implication $(b) \Rightarrow (a)$ follows from
        Seidenberg's Lemma, see Proposition 3.7.15 in \cite{kreuzer2000computational}.
    \end{proof}

    This proposition motivates the following definition.

    \begin{definition}
        A 0-dimensional radical ideal $I$ in $P$ is called \textbf{separable} if $I$ remains radical over every
        extension field of~$K$.
    \end{definition}

    Note that over a perfect field every radical ideal is separable. Given an ideal $I$ in $K[x_1, \dots, x_n]$,
    the field extensions $L$ of $K$ such that the radical of $IL[x_1, \dots, x_n]$ is separable can be characterized
    as follows.

    \begin{proposition}\label{prop:quasi_perfect}
        Let $R = P/I$ be a 0-dimensional affine $K$-algebra. For an extension field $L$ of $K$ the following
        are equivalent.
        \begin{enumerate}[(a)]
            \item The radical $J$ of $IL[x_1, \dots, x_n]$ is separable.
            \item For all $a \in R$ the minimal polynomial $\mu_a \in K[z]$ satisfies $\mysep(\mu_a) \in L[z]$.
        \end{enumerate}
    \end{proposition}
    \begin{proof}
        To prove (a) implies (b), let $a \in K[x_1, \dots, x_n]$. By Lemma~\ref{lemma:min_poly} the minimal polynomial
        $\mu_a \in K[z]$ of~$a+I$ and the minimal polynomial of $a+IL[x_1, \dots, x_n]$ coincide. Let $J$ be the radical
        of $IL[x_1, \dots, x_n]$, and let
        $$\varphi : L[x_1, \dots, x_n] \rightarrow L[x_1, \dots, x_n]/J$$
        be the canonical homomorphism. Let $h \in L[z]$ be the minimal polynomial of~$\varphi(a)$. Then $h$ has to
        divide $\mu_a$ since $0 = \varphi(\mu_a(a)) = \mu_a(\varphi(a))$.
        By Proposition~\ref{prop:sep_algebras} the polynomial $h$ has to be separable. Hence, $h$ divides
        $\mysep(\mu_a)$. Because $h(a) \in J$, we have $h^k(a) \in IL[x_1, \dots, x_n]$ for some $k>0$. This shows
        $\mysep(\mu_a) \mid h^k$ and therefore $\mysep(\mu_a) \mid h$ since $\mysep(\mu_a)$ is squarefree.

        To show the other implication, let $f_i \in I \cap K[x_i]$ for $i=1, \dots, n$. Then the ideal
        $\langle \mysep(f_i) \mid i =1, \dots, n \rangle + IL[x_1, \dots, x_n]$ is separable and the radical
        of $IL[x_1, \dots, x_n]$ by Seidenberg's Lemma (Proposition 3.7.15 in \cite{kreuzer2000computational}).
    \end{proof}

    \begin{definition}
        Let $R = P/I$ be a 0-dimensional affine $K$-algebra. An extension field $L$ of $K$ such that
        $L[x_1, \dots, x_n]/IL[x_1, \dots, x_n]$ satisfies the equivalent conditions in
        Proposition~\ref{prop:quasi_perfect} is called a \textbf{quasi-perfect} field for~$R$.
    \end{definition}

    A quasi-perfect field can be determined as follows.

    \begin{algorithm}{\textbf{(Computing a Quasi-Perfect Field)}}\\\label{alg:quasi-perfect}
        Let $R=P/I$ be a 0-dimensional affine $K$-algebra. The following sequence of instructions forms an algorithm
        which computes a quasi-perfect field for $R$.
        \begin{enumerate}[(1)]
            \item For $i=1, \dots, n$ compute the minimal polynomial $\mu_{x_i} \in K[z]$ of $\overline{x}_i$ in $R$.
            \item For $i=1, \dots, n$ compute the separable part $\mysep(\mu_{x_i}) \in L_i[z]$ where $L_i$ is an
                extension field of $K$.
            \item Determine a common extension field $L$ of $K$ such that $L_i \subseteq L$ for all
                $i=1, \dots, n$ and return it.
        \end{enumerate}
    \end{algorithm}

    \begin{proof}
        Let $J$ be the radical of $I L[x_1, \dots, x_n]$. Since $\mysep(\mu_{x_i}) \in L[z]$
        for $i=1, \dots, n$, the elements $\overline{x}_1, \dots, \overline{x}_n$ in $L[x_1, \dots,x_n]/J$
        are separable. Hence, $J$ is separable by Proposition~\ref{prop:sep_algebras}, proving that $L$
        is quasi-perfect for $R$.
    \end{proof}

    Independent of the base field the separable elements of a 0-di\-men\-sional affine $K$-algebra always form a
    subalgebra.

    \begin{proposition}\label{prop:separable_alg}
        Let $K$ be a field and $R$ a 0-dimensional affine $K$-algebra. The set $S$ of all separable elements
        of $R$ is a $K$-subalgebra of~$R$.
    \end{proposition}
    \begin{proof}
        \cite{kreuzer2016computational}, Proposition~5.5.3
    \end{proof}

    Note that Example~5.5.4 in~\cite{kreuzer2016computational} uses a definition of separability which is not
    applicable to non-perfect fields. It does therefore not provide a counterexample to
    Proposition~\ref{prop:separable_alg} in the case of non-perfect fields.

    \begin{definition}
        Let $R$ be a 0-dimensional affine $K$-algebra. The $K$-subalgebra of~$R$ which consists of all separable
        elements is called the \textbf{separable subalgebra} of~$R$ and is denoted by $R^{\mysep}$.
    \end{definition}

Over a perfect field~$K$, we always have a direct decomposition of~$R$ into its separable subalgebra 
and the zero radical. Over a non-perfect field it might be necessary to extend the base field 
to a quasi-perfect field for~$R$ to obtain such a decomposition.

    \begin{proposition}\label{prop:decomp}
        Let $R$ be a 0-dimensional affine $K$-algebra, and assume that $K$ is a quasi-perfect field for $R$.
        Every element $a\in R$ has a unique decomposition $a=b+r$ with $b \in R^{\mysep}$ and $r \in \Rad(0)$.
        In particular, we have a decomposition $R = R^{\mysep} \oplus \Rad(0)$ into a direct sum of $K$-vector
        subspaces.
    \end{proposition}
    \begin{proof}
        Let $\mu_a$ be the minimal polynomial of $a$. By Proposition~\ref{prop:quasi_perfect} we have
        $\mysep(\mu_a) \in K[x]$. Using this, the claim is a straightforward generalization of
        Proposition~5.5.6 in~\cite{kreuzer2016computational}.
    \end{proof}

    \begin{definition}
        Let $a \in R$. In the decomposition $a = b+r$ with $b \in R^{\mysep}$ and $r \in \Rad(0)$ the element
        $b$ is called the \textbf{separable part} of~$a$ and is denoted by $a^{\mysep}$. The element $r$ is called the
        \textbf{nilpotent part} of~$a$ and is denoted by~$a^{\nil}$.
    \end{definition}

    Proposition~5.5.6 in~\cite{kreuzer2016computational} yields an algorithm for computing this decomposition.

    \begin{algorithm}{\textbf{(Computing the Separable and the Nilpotent Part)}}\\\label{alg:sep_dec}
        Let $R$ be a 0-dimensional affine $K$-algebra and let $a \in R$. Assume that $K$ is a quasi-perfect field
        for $R$. Consider the following sequence of instructions.
        \begin{enumerate}[(1)]
            \item Compute the minimal polynomial $\mu_a \in K[x]$ of $a$, and compute $f = \mysep(\mu_a)$.
            \item Let $i=0$, $b_0 = a$ and $r_0 = 0$.
            \item Increase $i$ by one, let $b_i = b_{i-1} - \frac{f(b_{i-1})}{f'(b_{i-1})}$, and let
            $r_i = r_{i-1} - \frac{f(b_{i-1})}{f'(b_{i-1})}$.
            \item Repeat Step~(3) until $f(b_i) = 0$.
            \item Return the pair $(b_i, r_i)$.
        \end{enumerate}
        This is an algorithm which computes a pair $(b,r)$ such that $b$ is the separable part and $r$ is the nilpotent
        part of~$a$.
    \end{algorithm}

    We conclude the first part of this section by citing a useful result on the structure of 0-dimensional affine
    $K$-algebras. Once again the requirement that $K$ is a perfect field can be weakened to the assumption that $K$ is
    a quasi-perfect field for the given algebra.

    \begin{proposition}\label{prop:separable_decomp}
        Let $R$ be a 0-dimensional affine $K$-algebra such that $K$ is a quasi-perfect field for $R$.
        \begin{enumerate}[(a)]
            \item We have $R^{\mysep} \cong K[\overline{x}_1^{\mysep}, \dots, \overline{x}_n^{\mysep}]$
            \item The $K$-algebra homomorphism $R \rightarrow R^{\mysep}$ given by $a \mapsto a^{\mysep}$ induces a
                $K$-algebra isomorphism $R/\Rad(0) \cong R^{\mysep}$.
        \end{enumerate}
    \end{proposition}
    \begin{proof}
        \cite{kreuzer2016computational}, Proposition~5.5.11 and Corollary~5.5.13
    \end{proof}

%%%%%%%%%%%%%%%%%%%%%%%%%%%%%%%%%
% Subsection 3.2
%%%%%%%%%%%%%%%%%%%%%%%%%%%%%%%%%
\bigbreak
\subsection*{Computing Exponent Lattices in 0-Dimensional Algebras}

Let us now turn to the problem of computing the exponent lattice of units in a 0-dimensional
    affine $K$-algebra $R$. To be specific, for elements $f_1, \dots, f_n \in R^\times$ our goal is to compute a basis
    of the lattice given by all $a \in \mathbb{Z}^n$ with $f_1^{a_1} \cdots f_n^{a_n} = 1$.
    In 0-dimensional $\mathbb{Q}$-algebras, exponent lattices can be computed using the algorithm presented
    in Section~8 of \cite{lenstra2018algorithms}. Here we present a new algorithm for 0-dimensional algebras in
    finite characteristic. Additionally we present a slightly generalized version of the algorithm
    by Lenstra and Silverberg.

    \begin{proposition}\label{prop:isomorphism}
        Let $R$ be a 0-dimensional affine $K$-algebra, let $r \in R^\times$, and assume that $K$ is a quasi-perfect
        field for $R$.
        \begin{enumerate}[(a)]
            \item The set $1+\Rad(0) \coloneqq \{1+r \mid r \in \Rad(0)\}$ is a subgroup of $R^\times$.
            \item The map $$\varphi : R^{\times} \rightarrow (1+\Rad(0))\times (R^{\mysep})^\times, \quad r \mapsto
                (r\cdot (r^{\mysep})^{-1}, r^{\mysep})$$ is a group isomorphism.
            \item Let $m \in \mathbb{N}$ such that $(r^{\nil})^m = 0$. Then we have
                $r(r^{\mysep})^{-1} = \sum_{i=0}^{m-1} r^{-i}(r^{nil})^i$.
        \end{enumerate}
    \end{proposition}
    \begin{proof}
        Clearly, the elements in $1+\Rad(0)$ form a group with respect to multiplication. Since every
        element $r \in \Rad(0)$ is nilpotent, we have $(1-r)^{-1} = \sum_{i=1}^m r^i$ for $m$ large enough.
        Therefore $1+\Rad(0)$ is a subgroup of $R^\times$. This proves~(a).

        The map $\varphi$ in (b) is a well-defined group homomorphism since the nilpotent part and the separable part
        of an element $r$ are uniquely determined by Proposition~\ref{prop:decomp} and $R^{\mysep}$ is a subalgebra
        of~$R$. To show the surjectivity of $\varphi$, let $(1+r,s) \in (1+\Rad(0))\times (R^{\mysep})^\times$. Then
        $rs \in \Rad(0)$ and we have
        $$\varphi(rs+s) = ((rs+s)s^{-1}, s) = (r+1, s).$$
        Finally, let $r,r' \in \Rad(0)$ and $s,s' \in R^{\mysep}$. Then $\varphi(r+s) = \varphi(r'+s')$ implies
        $s = s'$ and $rs^{-1}+1 = r's{-1}+1$. Hence $r=r'$, and we proved the injectivity of $\varphi$.

        Part~(c) follows from $(r^{\mysep})^{-1} = (r-r^{\nil})^{-1} = \sum_{i=0}^{m-1} r^{-i-1}(r^{\nil})^{-1}$.
    \end{proof}

    The isomorphism in Part~(b) of the proposition allows us to compute the exponent lattices in $R^{\mysep}$ and in
    $1+\Rad(0)$ separately. Using the Chinese remainder theorem $R^{\mysep}$ can be further decomposed into a direct
    product of fields. Methods for computing exponent lattices in fields are given in
    Section~\ref{sec2}. It remains to solve the problem in~$1+\Rad(0)$. In characteristic zero
    we can use the following result.

    \begin{proposition}\label{prop:log}
        Let $K$ be a field of characteristic zero, and let $R$ be a 0-dimensional affine $K$-algebra. Assume that
        $m \in \mathbb{N}$ such that $\Rad(0)^m = 0$. Then the map
        $$\log : 1+ \Rad(0) \rightarrow \Rad(0), \quad 1+v \mapsto \ssum_{i=1}^{m-1} \frac{v^i}{i},$$
        is a group isomorphism from the multiplicative group $1+ \Rad(0)$ to the additive group of $\Rad(0)$.
    \end{proposition}
    \begin{proof}
        \cite{lenstra2018algorithms}, Proposition 8.1
    \end{proof}

    In characteristic $p$ such a bijective logarithm map can only be the trivial homomorphism. An inverse map
    $\myexp : \Rad(0) \rightarrow 1+\Rad(0)$ would have to satisfy $\myexp(0) = 1$, $\myexp(a+b) = \myexp(a)\myexp(b)$
    and therefore
    $$(\myexp(a)-1)^p = \myexp(a)^p -1 = \myexp(pa) -1 = \myexp(0)-1 = 0$$
    for all $a,b \in \Rad(0)$. Hence $\myexp$ can only be the trivial homomorphism. Instead we can use that in positive
    characteristic every element in $1+\Rad(0)$ has finite order.

    \begin{remark}\label{remark:group_structure}
        Let $K$ be a field of characteristic $p \geq 2$ and $R$ a 0-dimensional affine $K$-algebra. For
        $r \in \Rad(0)$ let $m$ be the nilpotency index of $r$, i.e., the smallest integer $m>0$ such that $r^m=0$.
        If $k>0$ is the smallest integer such that $p^k \geq m$, then the order of the element $1+r$ in the group
        $1+\Rad(0)$ is $p^k$. The problem of computing the exponent lattice of $(1+r_1, \dots, 1+r_s) \in (1+\Rad(0))^s$
        can therefore be restated as computing the group structure of the finite
        abelian $p$-group generated by $1+r_1, \dots, 1+r_s$. This can be achieved using one of the algorithms
        in~\cite{buchmann1997some,buchmann2005computing,teske1998space}.
    \end{remark}

    Since there are substantial differences depending on the characteristic of $K$, we present two separate algorithms
    for computing the exponent lattice. The next algorithm is a slightly generalized version of Algorithm~8.3
    from~\cite{lenstra2018algorithms}.

    \begin{algorithm}{\textbf{(Computing Exponent Lattices in Characteristic Zero)}}\\\label{alg:zero_dim_lattice0}
        Let $K$ be a field of characteristic zero and $R$ a 0-dimensional affine $K$-algebra.
        The following sequence of instructions forms an algorithm which computes the exponent lattice of
        $(f_1, \dots, f_k) \in (R^\times)^k$.
        \begin{enumerate}[(1)]
            \item Compute the maximal ideals $\mathfrak{m}_1, \dots \mathfrak{m}_s$ of $R$.
            \item For $i=1, \dots, s$ compute the exponent lattice~$\Lambda_i$ of
                $(\overline{f}_1, \dots, \overline{f}_k)$ where $\overline{f}_i$ is the canonical image of
                $f_i$ in the field $R/\mathfrak{m}_i$.
            \item Using Algorithm~\ref{alg:sep_dec} compute the decomposition $f_i = f_i^{\mysep}+f_i^{\nil}$
                for $i=1, \dots, k$.
            \item For $i=1, \dots, k$ compute $\log(f_i (f_i^{\mysep})^{-1}) =
                \log(\sum_{j=0}^{m_i-1} f_i^{-j}(f_i^{\nil})^j)$ where $m_i$ is the nilpotency index of $f_i^{\nil}$.
            \item Let $y_1, \dots, y_k$ be indeterminates, and consider the linear equation in the free abelian
                group $\Rad(0)$ given by
                $$y_1 \log(f_1 (f_1^{\mysep})^{-1}) + \cdots + y_k \log(f_k (f_k^{\mysep})^{-1}) = 0.$$
                After clearing denominators, this yields a homogeneous linear equation over $\mathbb{Z}$.
                Compute its solution space $M \subseteq \mathbb{Z}^k$.
            \item Compute the lattice $M \cap \Lambda_1 \cap \cdots \cap \Lambda_s$ and return it.
        \end{enumerate}
    \end{algorithm}
    \begin{proof}
        Let $a = (a_1, \dots, a_k) \in \mathbb{Z}^{k}$, and let
        $\varphi : R^{\times} \rightarrow (1+\Rad(0))\times (R^{\mysep})^\times$ be the isomorphism given as in
        Proposition~\ref{prop:isomorphism}. Then the tuple $a$ is an element of the exponent
        lattice of $(f_1, \dots, f_k)$ if and only if
        $$ \varphi(f_1)^{a_1} \cdots \varphi(f_k)^{a_k} = (f_1 (f_1^{\mysep})^{-1}, f_1^{\mysep})^{a_1} \cdots
        (f_k (f_k^{\mysep})^{-1}, f_k^{\mysep})^{a_k} = 1.
        $$
        Part (b) of Proposition~\ref{prop:separable_decomp} and the Chinese Remainder Theorem then imply
        $R^{\mysep} \cong R/\mathfrak{m}_1 \times \cdots \times R/\mathfrak{m}_s$.
        Therefore the product of the second components equals 1 if and only if $a \in \Lambda_i$
        for all $i=1 \dots, s$. Applying the logarithm map, we then get that $a$ is in the exponent lattice of
        $(f_1 (f_1^{\mysep})^{-1}, \dots, f_1 (f_1^{\mysep})^{-1})$ in $1+\Rad(0)$ if and only if $a$ is a solution of
        the linear system in Step~(5).
    \end{proof}

    The exponent lattices in Step~(3) of this algorithm can be computed as described in Section~\ref{sec2}.
    Let us see an example which illustrates this algorithm.

    \begin{example}
        Let $P = \mathbb{Q}(z)[x^\pm,y^\pm]$, and consider the 0-di\-men\-sional ideal
        $$I = \langle 3x -y -2z, y^2 -2zy +z^2 \rangle \subseteq P.$$
        Let $R=P/I$, and let $f_1,f_2$ be the residue classes of $x,y$ in $R$ and let $f_3=z$. We follow the steps of
        Algorithm~\ref{alg:zero_dim_lattice0} to compute the exponent lattice of $(f_1,f_2,f_3)$ in~$R$. The only maximal
        ideal of $R$ is $\mathfrak{m} = \langle \overline{y}-z, \overline{x}-z \rangle$. Since it is
        binomial, we easily obtain the exponent lattice $\Lambda = \langle (0,1,-1), (1,0,-1) \rangle$ of
        $(f_1, f_2, f_3)$ in $P/\mathfrak{m}$.
        Next, we compute $a^{\mysep} = b^{\mysep} = z$, $a^{\nil} = \overline{x}-z$ and
        $b^{\nil} = \overline{y}-z$. Then, after clearing denominators, we get
        \begin{align*}
            \log(f_1(f_1^{\mysep})^{-1}) &= -z\overline{y}^{-1} +1,\\
            \log(f_2(f_2^{\mysep})^{-1}) &= -3z\overline{y}^{-1} +3 \\
            \log(f_3(f_3^{\mysep})^{-1}) &= \log(1) = 0.
        \end{align*}
        This yields the homogeneous linear system of equations over~$\mathbb{Z}$ given by the matrix
        $$\left(\begin{array}{rrr}
                            1 & 3 & 0 \\
                            -1 & -3 & 0
                        \end{array}\right).$$
        Its kernel is given by $M = \langle (-3,1,0), (0,0,1) \rangle$. From this we then obtain the exponent lattice
        $\Lambda \cap M = \langle (3, -1, -2) \rangle$.
    \end{example}

    In finite characteristic we can compute the exponent lattice as follows.

    \begin{algorithm}{\textbf{(Computing Exponent Lattices in Finite Characteristic)}}\\\label{alg:zero-dimP}
        Let $K$ be a field of prime characteristic $p$ and $R$ a 0-di\-men\-sional affine $K$-algebra.
        The following sequence of instructions forms an algorithm which computes the exponent lattice of
        $(f_1, \dots, f_k) \in (R^\times)^k$.
        \begin{enumerate}[(1)]
            \item Using Algorithm~\ref{alg:quasi-perfect} compute a quasi-perfect field $L$ for $R$.
            \item Compute the maximal ideals $\mathfrak{m}_1, \dots \mathfrak{m}_s$ of
                $S = L \otimes_K R$.
            \item For $i=1, \dots, s$ compute the exponent lattice~$\Lambda_i \subseteq \mathbb{Z}^k$ of
                $(\overline{f}_1, \dots, \overline{f}_k)$ where $\overline{f}_j$ is the canonical image of $f_j$
                in the field~$S/\mathfrak{m}_i$.
            \item For $i = 1, \dots, k$ let $g_i$ be the canonical image of $f_i$ in $S$. Using
                Algorithm~\ref{alg:sep_dec} compute the decompositions $g_i = g_i^{\mysep}+g_i^{\nil}$.
            \item Let $h_i = 1+g_i^{\nil}(g_i^{\mysep})^{-1}$. Using Remark~\ref{remark:group_structure} compute the
                exponent lattice of $(h_1, \dots, h_k)$ in $1+\Rad(0)$ where $\Rad(0)$ is the zero radical of~$S$.
            \item Compute the lattice $M \cap \Lambda_1 \cap \cdots \cap \Lambda_s$ and return it.
        \end{enumerate}
    \end{algorithm}

    \begin{proof}
        Clearly, a tuple $a = (a_1,\dots,a_k) \in \mathbb{Z}^k$ is in the exponent lattice of $(f_1,\dots,f_k)$ if and
        only if $g_1^{a_1} \cdots g_k^{a_k} = 1$ in $S$. Let
        $\varphi : S^{\times} \rightarrow (1+\Rad(0))\times (S^{\mysep})^\times$ be the isomorphism given as in
        Proposition~\ref{prop:isomorphism}. Then $a$ is in the exponent lattice of
        $(g_1, \dots, g_k)$ if and only if
        $$
        \varphi(g_1)^{a_1} \cdots \varphi(g_k)^{a_k} = (g_1 (g_1^{\mysep})^{-1}, g_1^{\mysep})^{a_1} \cdots
        (g_k (g_k^{\mysep})^{-1}, g_k^{\mysep})^{a_k} = 1.
        $$
        Part (b) of Proposition~\ref{prop:separable_decomp} and the Chinese Remainder Theorem then imply
        $S^{\mysep} \cong S/\mathfrak{m}_1 \times \cdots \times S/\mathfrak{m}_s$.
        Therefore the product of the second components is equal to~1 if and only if $a \in \Lambda_i$ for all
        $i=1 \dots, s$. The product of the first components is equal to~1 if and only if $a \in M$. We therefore
        conclude that $a$ is an element of the exponent lattice of $(f_1, \dots, f_k)$ if and only if
        $a \in M \cap \Lambda_1 \cap \cdots \cap \Lambda_s$.
    \end{proof}

    Again, the exponent lattices in Step~(3) of this algorithm can be computed as described in Section~\ref{sec2}.

    \begin{example}
        Let $P = \mathbb{F}_5[x^\pm,y^\pm]$ and consider the ideal
        $$I =\langle -x+y-1, x^2-x-1 \rangle \subseteq P.$$ Let us compute
        the exponent lattice of $(f_1, f_2)$ where $f_1 = \overline{x}$ and $f_2 = \overline{y}$ in $R=P/I$. Extending
        the base field is not necessary since $\mathbb{F}_5$ is a perfect field. The ideal $I$ is primary and its
        radical is $\langle y+1, x+2 \rangle$. The nilpotent part of $f_1$ is $\overline{x}+2$ and its separable part is
        $-2$. Similarly, the nilpotent part of $f_2$ is $y+1$ and its separable part is $-1$. The exponent
        lattice of the separable parts in the field $\mathbb{F}_5[x^\pm, y^\pm]/\langle y+1, x+2 \rangle$ is given
        by $\Lambda = \langle (2,-1), (0,2) \rangle$. Next we compute
        \begin{align*}
            h_1 &= 1+f_1^{\nil}(f_1^{\mysep})^{-1} = -3x \\
            h_2 &= 1+f_2^{\nil}(f_2^{\mysep})^{-1} = -y
        \end{align*}
        The nilpotency index both for $h_1-1$ and for
        $h_2-1$ is 2, which means the order of $h_1$ and $h_2$ in $1+\Rad(0)$ is 5. Hence a tuple
        $(a_1, a_2) \in \mathbb{Z}^2$ such that $h_1^{a_1} \cdot h_2^{a_2} = 1$ has to satisfy $0 \leq a_1, a_2 \leq 5$.
        Exhaustive search then yields the exponent lattice $M = \langle (1,2), (0,5) \rangle$ of the elements $h_1, h_2$.
        Finally, we get $\Lambda \cap M = \langle (2,-1), (0,10) \rangle$ as the exponent lattice of $(f_1, f_2)$.
        These computations show that the unitary binomial part of $I$ is given by $\langle x^2-y, y^{10}-1 \rangle$.
    \end{example}

\bigskip\bigbreak
%%%%%%%%%%%%%%%%%%%%%%%%%%%%%%%%%%%%%%%%%%%%%%%%%%%%%%
%
%  Section 4: Unit Lattices and their Associated Characters
%
%%%%%%%%%%%%%%%%%%%%%%%%%%%%%%%%%%%%%%%%%%%%%%%%%%%%%%

\section{Unit Lattices and their Associated Characters}

In the following let $K$ be a field, $P = K[x_1, \dots, x_n]$ the polynomial ring over~$K$, and $I$ an ideal in~$P$.
In this section we study unit lattices in affine algebras $R=P/I$ and their associated characters.
Let $F = (f_1, \dots, f_k)$ be a tuple of elements in a ring. For a tuple $a = (a_1, \dots, a_k) \in \mathbb{N}^k$
we use the shorthand notation $F^a = f_1^{a_1} \cdots f_k^{a_k}$. Analogously we use this notation for tuples
$a \in \mathbb{Z}^k$ if the elements $f_i$ are invertible. Additionally, we let
$a=a^+ - a^- = (a^+_1 - a^-_1,\dots, a^+_k-a^-_k)$ be the unique decomposition with $a^+_i = \max\{a_i,0\}$ and
$a^-_i = \min\{a_i,0\}$.

\begin{proposition}\label{prop:unit_lattice}
    Let $F = (f_1, \dots, f_k)$ be a tuple of polynomials in $P$, let $I \subsetneq P$ be an ideal, and let
    $G$ be a subgroup of $K^\times$.
    \begin{enumerate}[(a)]
        \item If $I : \langle f_1 \cdots f_k \rangle^\infty = I$, then the residue classes of $f_1, \dots, f_k$ are
            non-zerodivisors in~$P/I$.
        \item The set
            $$
            \Lambda \coloneqq \left\{a \in \mathbb{Z}^k \mid F^{a^+} - g_a F^{a^-} \in I
            \text{ for some } g_a \in G\right\}
            $$
            is a lattice.
        \item If the residue classes of the elements $f_1, \dots, f_k$ are non-zerodivisors in $P/I$,
            then the map $\rho : \Lambda \rightarrow G$ given by $\rho(a) = g_{a}$
            for all~$a \in \Lambda$ is a well-defined group homomorphism.
    \end{enumerate}
\end{proposition}

\begin{proof}
    Let $g \in P$ such that $gf_i \in I$, then the assumption immediately implies $g \in I$. This proves~(a).

    To prove~(b) let $a, b \in \Lambda$. Then we have $F^{a^+}- g_a F^{a^-} \in I$ and $F^{b^+}-g_b F^{b^-} \in I$
    for some $g_a, g_b \in G$. Now $F^{a^-}-g_a^{-1} F^{a^+} \in I$ immediately implies $-a \in I$. We also have
    $F^{a^+}F^{b^+}-g_a g_b F^{a^-}F^{b^-} \in I$ and therefore $a+b \in \Lambda$.

    For Part~(c) assume that $F^{a^+}-gF^{a^-}$ and $F^{a^+}-g'F^{a^-}$ are in $I$. Then we have
    $-gF^{a^-}+g'F^{a^-} = (-g+g')F^{a^-} \in I$. Since $f_1, \dots, f_k$ are non-zerodivisors
    in $P/I$ this shows $g = g'$. Hence $g_a$ is uniquely determined and $\rho$ is a well-defined
    group homomorphism.
\end{proof}

In the following the condition in Part~(a) of this proposition, that $I$ is saturated with respect to the product
$f_1 \cdots f_k$ is crucial. This can be checked using one of the equivalent conditions given below.

\begin{remark}
    For an ideal $I \subseteq P$ the following are equivalent.
    \begin{enumerate}
        \item[(a)] $I : \langle f_1 \cdots f_k \rangle^\infty = I$
        \item[(b)] $I : \langle f_1 \cdots f_k \rangle = I$
        \item[(c)] $I : \langle f_i \rangle = I$ for $i=1, \dots, k$.
    \end{enumerate}
    This follows from the formulas
    $I : \langle f_1 \cdots f_k \rangle^\infty = \bigcup_{i \geq 1} (I : \langle f_1 \cdots f_k \rangle^i)$ and
    $$
    I : \langle f_1 \cdots f_k \rangle^\infty = (\cdots((I : \langle f_1 \rangle^\infty) :
    \langle f_2 \rangle^\infty) \cdots) : \langle f_K \rangle^\infty.
    $$
\end{remark}

Proposition~\ref{prop:unit_lattice} motivates the following definition.

\begin{definition}
    Let $F=(f_1,\dots,f_k) \in P^k$ be a tuple of polynomials, let $I \subseteq P$ be an ideal such that
    $I : \langle f_1 \cdots f_k \rangle = I$, and let $G$ be a subgroup of $K^\times$.

    \begin{enumerate}[(a)]
        \item The lattice consisting of all $a=(a_1,\dots,a_k) \in \mathbb{Z}^k$ such that
            $$
            f_1^{a^+_1} \cdots f_k^{a^+_k} \;-\; g_a\cdot f_1^{a^-_1} \cdots f_k^{a^-_k} = 0
            \quad\hbox{\rm for \ }g_a\in K^\times
            $$
            is called the \textbf{unit lattice of $F$ modulo $I$ with respect to $G$}.
        \item The group homomorphism $\rho : \Lambda \rightarrow K^\times$ given by $\rho(a) = g_a$
            is called its \textbf{associated character}.
    \end{enumerate}
\end{definition}

If the group $G$ in this definition is the whole group of units $K^\times$ then we will simply refer to this
lattice as the unit lattice of $F$ modulo $I$.

\begin{remark}\label{remark:character_data}
    In the following our goal is to compute the unit lattice~$\Lambda$ and the associated character~$\rho$ of a tuple
    $(f_1, \dots, f_k)$ modulo $I$. By this we mean computing a basis $b_1, \dots, b_m$ of $\Lambda$ together with
    elements $g_1, \dots, g_m \in K^\times$ such that $\Lambda(b_i) = g_i$.
    From now on we assume that unit lattices and their associated characters are given in this form.
\end{remark}

When computing the unit lattice of a tuple $(f_1, \dots, f_k)$ modulo $I$ it turns out to be useful to write $I$
as an intersection of ideals. It is then necessary to compute the following.

\begin{definition}
    Let $\Lambda,M \subseteq \mathbb{Z}^k$ be lattices, and let $K$ be a field. For characters
    $\rho : \Lambda \rightarrow K^\times$ and $\tau : M \rightarrow K^\times$ we call the lattice
    $$\{v \in \Lambda \cap M \mid \rho(v) = \tau(v)\}$$
    the \textbf{intersection of $\Lambda$ and $M$ with respect to $\rho$ and $\tau$}. We denote
    it by $(\Lambda, \rho) \cap (M, \tau)$.
\end{definition}

The intersection of lattices $\Lambda$ and $M$ in $\mathbb{Z}^k$ with respect to characters
$\rho : \Lambda \rightarrow K^\times$ and $\tau : M \rightarrow K^\times$ can be determined by computing the basis of
an exponent lattice in~$K^\times$.

\begin{algorithm}{\textbf{(Computing Lattice Intersections w.r.t. Characters)}}\\\label{alg:lattice_intersection}
    Let $\Lambda$, $M$, $\rho$ and $\tau$ be given as above. The following instructions form an algorithm which
    computes the intersection of $\Lambda$ and $M$ with respect to $\rho$ and $\tau$.
    \begin{enumerate}[(1)]
        \item Compute a basis $b_1, \dots, b_r \in \mathbb{Z}^k$ of the lattice $\Lambda \cap M$.
        \item Compute a basis $c_1, \dots, c_s \in \mathbb{Z}^r$ of the lattice of exponents
            $\mathcal{L} \subseteq \mathbb{Z}^r$ of
            $$\rho(b_1)\tau(b_1)^{-1}, \dots, \rho(b_r)\tau(b_r)^{-1} \text{ in } K^\times.$$
        \item For $i =1, \dots, s$ let $d_i = c_{i1}b_1 + \cdots + c_{ir}b_r$, and return the lattice
            $N = \langle d_1, \dots, d_s \rangle \subseteq \mathbb{Z}^k$.
    \end{enumerate}
\end{algorithm}
\begin{proof}
    Let $h_1, \dots, h_r \in \mathbb{Z}$. An element $a = h_1b_1+ \cdots + h_rb_r \in \Lambda \cap M$ satisfies
    $\rho(a) = \tau(a)$ if and only if
    \[\tau(b_1)^{h_1} \cdots \tau(b_r)^{h_r} = \rho(b_1)^{h_1} \cdots \rho(b_r)^{h_r}\]
    which is equivalent to $(h_1, \dots, h_r) \in \mathcal{L}$.
\end{proof}

A lattice $\Lambda$ in~$\mathbb{Z}^n$ together with a character $\rho:\; \Lambda \longrightarrow
K^\times$ yields a binomial ideal.

\begin{definition}
    Let $\Lambda \subseteq \mathbb{Z}^n$ be a lattice and $\rho:\; \Lambda \longrightarrow K^\times$ a character.
    The ideal
    $$
    I_{\Lambda,\rho} \;=\; \langle x_1^{a^+_1} \cdots x_n^{a^+_n} \;-\; \rho(a) \cdot x_1^{a^-_1} \cdots
    x_n^{a^-_n} \;\mid\; a = a^+ - a^- \in \Lambda \rangle
    $$
    in~$P$ is called the \textbf{lattice ideal} associated to $(\Lambda,\rho)$.
\end{definition}

For a detailed discussion of lattice ideals, see~\cite{villarreal2001monomial} or~\cite{herzog2018binomial}. In
the following we let $X = \{x_1, \dots x_n\}$ be the set of all indeterminates in $P$.
Proposition~\ref{prop:unit_lattice} in particular states that a binomial ideal $I$ that satisfies
$I : \langle x_1 \cdots x_n \rangle^\infty = I$ is a lattice ideal. The converse is also true.

\begin{proposition}\label{prop:lattice_to_binomial}
    \begin{enumerate}[(a)]
        \item A binomial ideal $I$ in $P$ is a lattice ideal if and only if it satisfies
            $I : \langle x_1 \cdots x_n \rangle^\infty = I$.
        \item Let $\Lambda \subseteq \mathbb{Z}^n$ be a lattice and let $\rho : \Lambda \rightarrow K^\times$ be a
            character. The lattice $\Lambda$ is generated by $b_1, \dots, b_k \in \mathbb{Z}^n$ if and only if
            $$I_{\Lambda, \rho} = \left\langle X^{b_i^+} - \rho(b_i) X^{b_i^-} \mid i=1, \dots, k \right\rangle :
            \langle x_1 \cdots x_n \rangle^\infty.$$
    \end{enumerate}
\end{proposition}

\begin{proof}
    \cite{villarreal2001monomial}, Theorem~8.2.8 and Lemma~8.2.11
\end{proof}

Given a lattice $\Lambda$ and a character $\rho : \Lambda \rightarrow K^\times$ as in Remark~\ref{remark:character_data},
Part~(b) of this proposition allows us to compute the corresponding lattice ideal. An alternative method which does not
use saturation is described in~\cite{hemmecke2009computing}.

\begin{corollary}\label{cor:lattice_ideal}
    Let $I \subseteq P$ be an ideal which satisfies $I : \langle x_1 \cdots x_n \rangle^\infty = I$. Let
    $\Lambda \subseteq \mathbb{Z}^n$ be the unit lattice of $(x_1, \dots, x_n)$
    modulo $I$, and let $\rho : \Lambda \rightarrow K^\times$ be its associated character. Then
    $\Bin(I) = I_{\Lambda, \rho}$. In particular $\Bin(I)$ is a lattice ideal.
\end{corollary}

\begin{proof}
    If $J \subseteq I$ is a binomial ideal, then the ideal $J: \langle x_1 \cdots x_n \rangle^\infty$ is contained
    in $I:\langle x_1 \cdots x_n \rangle^\infty$ and therefore also binomial. Hence
    $\Bin(I) = \Bin(I): \langle x_1 \cdots x_n \rangle^\infty$ is a lattice ideal by the proposition.
\end{proof}

    Let $f_1, \dots, f_k \in K[X]$ and $I$ an ideal in $K[X]$ with $I : \langle f_1 \cdots f_k \rangle = I$. The next
    result allows us to assume that the elements $f_1, \dots, f_k$ are invertible when computing their unit lattice and
    associated character.

    \begin{proposition}\label{prop:invertible}
        Let $F = (f_1, \dots, f_k)$ and $I \subseteq K[X]$ be given as above. Let $K[X]_F$ be the localization of
        $K[X]$ with respect to the multiplicatively closed set generated by the elements in~$F$. Then the unit lattices
        and associated characters of $F$ modulo $I$ and of $F$ modulo $IK[X]_F$ coincide.
    \end{proposition}
    \begin{proof}
        Proposition~\ref{prop:unit_lattice} yields $I : \langle f_1 \cdots f_k \rangle = IK[X]_F \cap K[X] = I$.
    \end{proof}

    When computing the unit lattice and the associated character modulo an ideal, it is sometimes necessary to extend
    the base field. This is possible by the following lemma. It is a straightforward generalization of Lemma~7
    in~\cite{jensen2017finding}. We provide a proof for the convenience of the readers.

    \begin{proposition}\label{prop:field_extension}
        Let $F = (f_1, \dots, f_k)$ and $I \subseteq K[X]$ be given as above, and let $L$ be an extension
        field of~$K$. Let $\Lambda$ be the unit lattice with associated character $\rho$ of $F$ modulo
        $I$, and let $\Lambda'$ be the unit lattice with associated character $\rho'$ of $F$ modulo
        $I L[X]$. Then we have $\Lambda = \Lambda'$ and $\rho(a) = \rho'(a) \in K$ for all $a \in \Lambda$.
    \end{proposition}
    \begin{proof}
        Clearly we have $\Lambda \subseteq \Lambda'$. Now assume $F^{a^+}-g F^{a^-} \in I L[X]$ for
        $a\in \mathbb{Z}^k$ and $g \in L^\times$. Since the residue classes of $f_1, \dots, f_k$ are non-zerodivisors in
        $K[X]/I$, their residue classes in $L[X]/I L[X]$ are also non-zerodivisors. Hence
        $\rho' : \Lambda' \rightarrow L^\times$ is a well-defined character, and we have
        $\rho'(a) = g$. To prove $g \in K$, notice that the assumption yields an expression
        $$
        F^{a^+}- g F^{a^-} = \ssum_i g_i F^{c_i} p_i
        $$
        with $a, c_i \in \mathbb{Z}^k$, $g_i \in \mathbb{K}$ and $p_i \in I$. For fixed $a,c_i$ and $p_i$
        this can be interpreted as a system of linear equations in the indeterminates~$g$ and~$g_i$. The
        coefficients of this system are in $K$ and it has a solution in $L$. Therefore the system also has a
        solution in $K$. This shows $F^{a^+}-g F^{a^-} \in I$, and since $g$ is uniquely determined we get
        $\rho'(a) = \rho(a) = g \in K$.
    \end{proof}

\bigskip\bigbreak
%%%%%%%%%%%%%%%%%%%%%%%%%%%%%%%%%%%%%%%%%%%%%%%%%
%
%  Section 5: Computing Unit Lattices
%
%%%%%%%%%%%%%%%%%%%%%%%%%%%%%%%%%%%%%%%%%%%%%%%%%

\section{Computing Unit Lattices}\label{section:computing_unit_lattices}

    Let $K$ be a field, $X = \{x_1, \dots, x_n\}$ a set of indeterminates, and $K[X]$ the polynomial ring over $K$ in
    the indeterminates $X$. Let $f_1, \dots, f_k \in K[X]$, and let $I \subseteq K[X]$ be an ideal such that the residue
    classes of $f_1, \dots, f_k$ in $K[X]/I$ are non-zerodivisors. This section is concerned with providing an
    algorithm for computing the unit lattice and the associated character of $f_1, \dots, f_k$ modulo~$I$.
    In the final section of this paper we then show how computing the binomial part of a general polynomial ideal
    reduces to computing unit lattices.

    The main idea is to reduce the computation to 0-dimensional ideals.
    Recall that a subset $U \subseteq X$ is said to be an independent set of indeterminates modulo~$I$ if we have
    $I \cap K[U] = \langle 0 \rangle$.

    \begin{proposition}\label{prop:indep_set}
        Let $I \subseteq K[X]$ be an ideal and $U \subseteq X = \{x_1, \dots, x_n\}$ be a maximal independent set of
        indeterminates modulo $I$. Let $\sigma$ be an elimination term ordering with respect to $X \setminus U$,
        and let $G$ be a Gr\"obner basis with respect to $\sigma$.
        \begin{enumerate}[(a)]
            \item The ideal $IK(U)[X \setminus U]$ is a 0-dimensional ideal.
            \item The set $G$ is a Gr\"obner basis of $IK(U)[X \setminus U]$.
            \item We have $IK(U)[X \setminus U] \cap K[X] = I : h^\infty$ where $h = \lcm\{\LC(g) \mid g \in G\}$ and $G$ is
                considered as a subset of $K(U)[X \setminus U]$.
        \end{enumerate}
    \end{proposition}
    \begin{proof}
        \cite{greuel2008singular}, Proposition~4.3.1
    \end{proof}

    This proposition together with the next lemma allows us to reduce the problem to computing
    unit lattices and associated characters modulo 0-dimensional ideals of the form~$IK(U)[X \setminus U]$.

    \begin{lemma}\label{lemma:splitting}
        Let $I \subseteq K[X]$ be an ideal and $f \in K[X]$ with $I : f^\infty = I : f^m$ for $m > 0$. Then we have
        \[I = (I : f^m) \cap \langle I, f^m \rangle.\]
    \end{lemma}
    \begin{proof}
        \cite{greuel2008singular}, Lemma~3.3.6
    \end{proof}

    Let us now present an algorithm for computing unit lattices and their associated characters modulo an ideal.

    \begin{algorithm}{\textbf{(Computing Unit Lattices and Associated Characters)}}\\\label{alg:unit_lattice}
        Let $I \subseteq K[X]$ be an ideal, and let $f_1, \dots, f_k \in K[X]$ such that
        $I : \langle f \rangle^\infty = I$ for $f=f_1 \cdots f_k$. Consider the following sequence of instructions.
        \begin{enumerate}[(1)]
            \item Compute a maximal independent set of indeterminates $U \subseteq X$ modulo $I$.
            \item Using Algorithm~\ref{alg:partial_char_zero_dim} compute the unit lattice $\Lambda$ and the associated
                character $\rho : \Lambda \rightarrow K^\times$ of
                $(f_1, \dots, f_k)$ modulo $IK(U)[X \setminus U] \cap K[X]$.
            \item Compute a Gr\"obner basis $G$ of $I$ with respect to an elimination ordering for the indeterminates
                in~$X \setminus U$.
            \item Compute $h = \lcm\{\LC(g) \mid g \in G\}$ with $G$ considered as a subset of $K(U)[X \setminus U]$,
                and compute $m > 0$ such that $I : h^\infty = I : h^m$.
            \item If $I : h^\infty \subseteq \langle I , h^m \rangle : \langle f \rangle^\infty$,
                return the unit lattice $\Lambda$ and the associated character $\rho$. Otherwise apply the
                algorithm recursively to $\langle I , h^m \rangle : \langle f \rangle^\infty$ and obtain a unit
                lattice $M$ with associated character~$\tau$.
            \item Apply Algorithm~\ref{alg:lattice_intersection} to compute the lattice
                $N = (\Lambda, \rho) \cap (M, \tau)$. Return $N$ together with the character~$\rho_{\mid N}$.
        \end{enumerate}
        This is an algorithm which computes the unit lattice and the associated character of $(f_1, \dots, f_k)$
        modulo~$I$.
    \end{algorithm}

    \begin{proof}
        We first show that the algorithm terminates. Since
        $I$ is contained in $J = (I+\langle h \rangle):\langle f \rangle^\infty$,
        the dimension of $J$ is less than or equal to the dimension of $I$. Furthermore the maximal
        independent set $U$ modulo $I$ is not independent modulo $J$ since $h \in K[U] \cap J$. This means in each
        recursive call of the algorithm the number of possible maximal independent sets or the dimension decreases.

        By Lemma~\ref{lemma:splitting} we have $I = (I : h^m) \cap \langle I, h^m \rangle$. Let us show that this
        implies $I = (I : h^m) \cap (\langle I, h^m \rangle : \langle f \rangle^\infty)$. An element
        $g \in (I : h^m) \cap (\langle I, h^m \rangle : \langle f \rangle^\infty)$ satisfies
        $gf^k \in \langle I, h^m \rangle$ for some $k \in \mathbb{N}$. Hence we get
        $gf^k \in (I : h^m) \cap (\langle I, h^m \rangle) = I$. This shows $g \in I$
        since we have $I : \langle f \rangle^\infty = I$ by assumption. Proposition~\ref{prop:indep_set} shows that the
        unit lattice and the associated character of $(f_1, \dots, f_k)$ modulo $I : h^m$ are given by $\Lambda$ and
        $\rho$. Therefore we conclude that the unit lattice of $(f_1, \dots, f_k)$ modulo $I$ is given by $N$, and
        that its associated character is $\rho_{\mid N}$.
    \end{proof}

    This algorithm already yields a method for computing the binomial part of ideals which are saturated with respect
    to the product of all indeterminates.

    \begin{corollary}{\textbf{(Computing the Binomial Part of Saturated Ideals)}}\\\label{cor:binomial_part_saturated}
        Let $I \subseteq K[X]$ be an ideal which satisfies $I : \langle x_1 \cdots x_n \rangle = I$. Then the following
        instructions form an algorithm which computes $\Bin(I)$.
        \begin{enumerate}[(a)]
            \item Using Algorithm~\ref{alg:unit_lattice} compute a basis $b_1, \dots, b_m \subseteq \mathbb{Z}^n$ of
                the unit lattice $\Lambda$ of $(x_1, \dots, x_n)$ modulo $I$ together with $c_1, \dots, c_m \in K^\times$
                such that the associated character $\rho : \Lambda \rightarrow K^\times$ satisfies $\rho(b_i) = c_i$
                for $i=1, \dots, m$.
            \item Return the ideal
                $$I_{\Lambda, \rho} = \left\langle X^{b_i^+} - c_i X^{b_i^-} \mid i=1, \dots, m \right\rangle :
            \langle x_1 \cdots x_n \rangle^\infty.$$
        \end{enumerate}
    \end{corollary}

    \begin{proof}
        By Corollary~\ref{cor:lattice_ideal} we have $\Bin(I) = I_{\Lambda, \rho}$. It then follows from
        Proposition~\ref{prop:lattice_to_binomial} that Step~(2) correctly computes the lattice ideal
        $I_{\Lambda, \rho}$.
    \end{proof}

    Later we will see that, when computing the binomial part of an arbitrary polynomial ideal, it is convenient to apply
    Algorithm~\ref{alg:unit_lattice} to a localized polynomial ring.

    \begin{remark}\label{remark:localization}
        Let $Y \subseteq X$ be a subset of indeterminates, $f_1, \dots, f_k \in K[X]$, and $I$ an ideal in $K[X]$.
        Assume that $I$ is saturated with respect to the product of all indeterminates in $Y$ and $f_1 \cdots f_k$.
        Then we have $IK[X]_Y \cap K[X] = I$ which means that the unit lattice and the associated character
        of $(f_1, \dots, f_k)$ modulo $I$ and modulo $IK[X]_Y$ coincide. The computations in $K[X]_Y$ can be performed
        in a polynomial ring by introducing a new indeterminate $z$ and adding the polynomial
        $x_{i_1} \cdots x_{i_m}z-1$ to $I$, where $Y = \{x_{i_1}, \dots, x_{i_m}\}$.
    \end{remark}

    Let us now investigate how in Step~(2) of Algorithm~\ref{alg:unit_lattice} the unit lattice and the associated
    character of $(f_1, \dots, f_k)$ modulo $IK(U)[X \setminus U] \cap K[X]$ can be computed. In a first step we show
    how the unit lattice with respect to $K(U)^\times$ modulo the zero-dimensional ideal $IK(U)[X \setminus U]$ can be
    found. In a second step we then show how from this the unit lattice of $(f_1, \dots, f_k)$ with respect to
    $K^\times$ modulo $IK(U)[X \setminus U] \cap K[X]$ can be obtained.

    Consider a 0-dimensional $K$-algebra $R$ and elements $f_1, \dots, f_k \in R^\times$. Our first goal is to
    compute the unit lattice together with the associated character of $(f_1, \dots, f_k)$ in~$R$. Using the following
    Proposition we can reduce the problem to computing exponent lattices in $R$. Depending on the characteristic
    of~$K$ these exponent lattices can then be computed using Algorithm~\ref{alg:zero_dim_lattice0} or
    Algorithm~\ref{alg:zero-dimP}. For $\mathbb{Q}$-algebras such a construction is given in
    Proposition~19 in~\cite{jensen2017finding}. But note that the construction in~\cite{jensen2017finding} is
    incorrect if the vector space dimension of the $\mathbb{Q}$-algebra is even.

\begin{proposition}\label{prop:reduction_to_unitary}
    Let $R$ be a 0-dimensional $K$-algebra, and let $f_1, \dots, f_k \in R^\times$. For $i=1, \dots, k$ consider
    the linear endomorphisms $\varphi_i : R \rightarrow R$ given by multiplication with $f_i$. Let $\ell$ be the
    dimension of $R$ as a vector space over $K$, and let $\zeta$ be a generator of the cyclic group of $\ell$-th roots
    of unity contained in $K$. We define $L$ to be the finite extension of $K$ which is obtained by adjoining all
    $\ell$-th roots of the determinants of the $\varphi_i$ to $K$. Let
    $$
    q_i = f_i/\sqrt[\ell]{\det(\varphi_i)} \in L \otimes_K R.
    $$
    Assume that $\Lambda \subseteq \mathbb{Z}^k$ is the unit
    lattice of $(f_1, \dots, f_k)$ with associated character $\rho : \Lambda \rightarrow K^\times$. Then
    a tuple $a \in \mathbb{Z}^k$ is in $\Lambda$ if and only if $a$ is in the exponent lattice of
    $(q_1, \dots, q_k, \zeta)$ projected onto the first $k$ components.
    In this case we have $\rho(a) = \zeta^b \prod_{i=1}^k \sqrt[\ell]{\det(\varphi_i)}^{a_i}$ for some $b \in \mathbb{Z}$.
\end{proposition}

\begin{proof}
    Let $a = (a_1, \dots, a_k) \in \Lambda$. Then we have $f_1^{a_1} \cdots f_k^{a_k} = g$ for some $g \in K^\times$,
    and the endomorphisms satisfy $\varphi_1^{a_1} \cdots \varphi_k^{a_k} = g \Id_R$.
    Taking determinants on both sides we get $\det(\varphi_1)^{a_1} \cdots \det(\varphi_k)^{a_k} = g^{\ell}$. Then
    taking $\ell$-th roots on both sides, we get $\prod_{i=1}^k \sqrt[\ell]{\det(\varphi_i)}^{a_i} = \zeta^b g$ for some
    $b \in \mathbb{Z}$. This implies
    $$
    q_1^{a_1} \cdots q_k^{a_k} \zeta^{-b} = \frac{1}{g}(f_1^{a_1} \cdots f_k^{a_k}) = 1.
    $$
    Conversely, let $q_1^{a_1} \cdots q_k^{a_k} \zeta^b =1$ and choose $g = \zeta^{-b}\prod_{i=1}^k \sqrt[\ell]{\det(\varphi_i)}^{a_i}$.
    Then we have $\prod_{i=1}^k f_i^{a_i} = g$.
\end{proof}

    When applying the reduction to 0-dimensional ideals using Algorithm~\ref{alg:unit_lattice}, we obtain
    0-dimensional ideals of the form $IK(U)[X \setminus U]$ where $U \subseteq X$ is a subset of indeterminates
    and $I$ is an ideal in $K[X]$. The unit lattice of $(f_1, \dots, f_k)$ modulo $IK(U)[X \setminus U]$ can then be
    computed using Proposition~\ref{prop:reduction_to_unitary}. However, the algorithm requires us to compute the unit
    lattice modulo $IK(U)[X \setminus U] \cap K[X]$ with respect to $K^\times$. In other words, we are not
    interested in all relations of the form $F^{a^+} - g F^{a^-}$ with $g \in K(U)^\times$ and $a \in \mathbb{Z}^k$ but
    only in those where $g$ is an element of~$K^\times$.

    \begin{algorithm}{\textbf{(Computing Unit Lattices in Zero-Dimensional Algebras)}}\label{alg:partial_char_zero_dim}
        Let $I \subseteq K[X]$ be an ideal and $U \subseteq X$ a subset of indeterminates such that
        $IK(U)[X \setminus U]$ is 0-dimensional. Let $F = \{f_1, \dots, f_k\}$ be a set of polynomials in $K[X]$
        such that $I : \langle f_1 \cdots f_k \rangle = I$. The following instructions form an algorithm which computes
        the unit lattice and the associated character of $(f_1, \dots, f_k)$ modulo the ideal
        $IK(U)[X \setminus U] \cap K[X]$ with respect to $K^\times$.
        \begin{enumerate}[(1)]
            \item Form the ring $R \coloneqq K(U)[X \setminus U]_F/I K(U)[X \setminus U]_F$, and for $i=1, \dots, k$
                let $\varphi_i$ be the linear endomorphisms of $R$ given by the multiplication with $f_i$.
            \item Let $\ell$ be the dimension of $R$ as a vector space over $K(U)$. Construct a finite extension
                $L$ of $K(U)$ which contains all the $\ell$-th roots of the determinants of all~$\varphi_i$.
            \item Determine a generator $\zeta$ of the cyclic group of $\ell$-th roots of unity contained in~$K$.
            \item For $i=1, \dots, k$ let $q_i = \overline{f}_i/\sqrt[\ell]{\det(\varphi_i)}$ where $\overline{f}_i$ is
                the residue class of $f_i$ in $L \otimes_{K(U)} R$, and compute the exponent lattice
                $M' \subseteq \mathbb{Z}^k$ of $(q_1, \dots, q_k, \zeta)$ in the 0-dimensional $\mathbb{K}$-algebra
                $L \otimes_{K(U)} R$.
            \item Compute the projection $M \subseteq \mathbb{Z}^k$ of $M'$ onto its first $k$-components.
            \item For $i =1, \dots, k$ write $\det(\varphi_i) = g_i p_i$ where $p_i \in K[U]$ is a monic polynomial
                and $g_i \in K^\times$.
            \item Using Algorithm~\ref{alg:exp_lattice_function_field} compute the exponent lattice
                $N \subseteq \mathbb{Z}^k$ of $(p_1^\ell, \dots, p_k^\ell)$ in~$K(U)^k$.
            \item Compute a basis $b_1, \dots, b_r \in \mathbb{Z}^k$ of $\Lambda = M \cap N$. For each $b_i$ compute the
                normal form $c_i \in K^\times$ of $f_1^{b_{i1}} \cdots f_k^{b_{ik}}$ modulo $I K(U)[X \setminus U]_F$.
            \item Return the lattice $\Lambda$ and the associated character $\rho: \Lambda \rightarrow K^\times$
                given by $b_i \mapsto c_i$ for $i=1, \dots, r$.
        \end{enumerate}
    \end{algorithm}
    \begin{proof}
        First, we note that by Proposition~\ref{prop:invertible} the unit lattice is stable under localization. It is
        therefore enough to compute the unit of lattice $(\overline{f}_1, \dots, \overline{f}_k)$ modulo
        the ideal $I K(U)[X \setminus U]_F$.
        Now let $a = (a_1, \dots, a_k) \in \mathbb{Z}^k$, and assume that $f_1^{a_1} \cdots f_k^{a_k} = g$ for some
        $g \in K^\times$. Then Proposition~\ref{prop:reduction_to_unitary} implies
        $q_1^{a_1} \cdots q_k^{a_k} \zeta^b = 1$ for some $b \in \mathbb{Z}$ and
        $\prod_{i=1}^k \sqrt[\ell]{\det(\varphi_i)}^{a_i} = \zeta^b g$. This shows $a \in M$ and
        $$(g_1 p_1)^{a_1} \cdots (g_k p_k)^{a_k} = g^\ell.$$
        Taking the leading coefficient of both sides of the equation we get $p_1^{a_1} \cdots p_k^{a_k} = 1$ and
        ${g_1}^{a_1} \cdots {g_k}^{a_k}= g^\ell$. We therefore conclude $a \in \Lambda = M \cap N$ and
        $g = \zeta^{-b} \prod_{i=1}^k \sqrt[\ell]{g_i}^{a_i}$.

        Conversely, let $a = (a_1, \dots, a_k) \in \Lambda$. Then for some $b \in \mathbb{Z}$ we have
        $$f_1^{a_1} \cdots f_k^{a_k} = \zeta^b \sprod_{i=1}^k \sqrt[\ell]{\det(\varphi_i)}^{a_i}\in K(U)$$
         by Proposition~\ref{prop:reduction_to_unitary}. Since $a \in N$ we also have
        $\prod_{i=1}^k \det(\varphi_i)^{a_i} = \prod_{i=1}^k {g_i}^{a_i}$. Taking the $\ell$-th root on both sides
        we get $f_1^{a_1} \cdots f_k^{a_k} = \zeta^b \prod_{i=1}^k \sqrt[\ell]{g_i}^{a_i} \in K'$, where $K'$ is a finite
        extension of $K$ which contains the $\ell$-th roots of the~$g_i$. Proposition~\ref{prop:field_extension} then
        shows $\zeta^b \prod_{i=1}^k \sqrt[\ell]{g_i}^{a_i} \in K^\times$. The normal forms $c_1, \dots, c_r$ computed
        in Step~(8) are therefore indeed elements of $K^\times$.
    \end{proof}

    If there exists an algorithm for factoring polynomials in $K[x]$, then an $\ell$-th root of unity as required in
    Step~(3) can be obtained by factoring $x^\ell-1 \in K[x]$. The linear factors then correspond to the $\ell$-th roots
    of unity contained in $K$.

\bigskip\bigbreak
%%%%%%%%%%%%%%%%%%%%%%%%%%%%%%%%%%%%%%%%%%%%%%%%%%%
%
%  Section 6: Computing the Binomial Part
%
%%%%%%%%%%%%%%%%%%%%%%%%%%%%%%%%%%%%%%%%%%%%%%%%%%%

\section{Computing the Binomial Part}

In this section we show how the computation of the binomial part of an ideal~$I$ in~$P$ can be reduced to computing
unit lattices and their associated characters. If an ideal $I$ satisfies $I : \langle x_1 \cdots x_n \rangle^\infty = I$,
then its binomial part can be computed as described in Corollary~\ref{cor:binomial_part_saturated}. If $I$ does not
satisfy this property, then our first step is to decompose it as follows.

\begin{definition}
Let $Y \subseteq \{x_1, \dots, x_n\}$. An ideal $I \subseteq P$ is called \textbf{$Y$-cellular} if the following 
conditions hold.
\begin{enumerate}
\item[(a)] $I = I : \langle \prod_{x_i \in Y} x_i \rangle^\infty$.

\item[(b)] For every $x_i \notin Y$, there exists an integer $d_i>0$ such that $x_i^{d_i} \in I$.
\end{enumerate}
\end{definition}

Every ideal can be written as an intersection of cellular ideals. This decomposition was first introduced
in~\cite{eisenbud1996binomial}. Using Lemma~\ref{lemma:splitting}, we obtain a straightforward algorithm for
decomposing an ideal into cellular parts (see also~\cite{castilla2000cellular}).

\begin{algorithm}{\textbf{(Computing Cellular Decompositions)}}\\\label{alg:cellular_decomposition}
Let $I \subseteq P$ be an ideal. Consider the following sequence of instructions.
\begin{enumerate}
\item[(1)] Determine an indeterminate $x_i$ such that $I : \langle x_i \rangle^\infty \neq \langle 1 \rangle$
and $I : \langle x_i \rangle^\infty \neq I$. If no such indeterminate exists, return $\{I\}$.

\item[(2)] Compute an integer $m>0$ such that $I : \langle x_i \rangle ^\infty = I : \langle x_i \rangle^m$.

\item[(3)] Recursively apply the algorithm to the ideals $I : \langle x_i \rangle^m$ and $\langle I, x_i^m \rangle$,
    and return the union of their cellular decompositions.
\end{enumerate}
This is an algorithm which computes a set of cellular ideals such that~$I$ is the intersection of the ideals
contained in the set.
\end{algorithm}

\begin{proof}
    The correctness of this algorithm follows from the formula
    $$
    I : \langle x_1 \cdots x_n \rangle^\infty = (\cdots((I : \langle x_1 \rangle^\infty) : \langle x_2 \rangle^\infty)
    \cdots) : \langle x_n \rangle^\infty.
    $$
\end{proof}

By definition, all monomials contained in an ideal~$I$ are also contained in $\Bin(I)$.

\begin{definition}
Let $I \subseteq P$ be an ideal. The ideal $\Mon(I)$ generated by all monomials contained in~$I$ is called
the \textbf{monomial part} of~$I$.
\end{definition}

Computing generators of the monomial part is decisively easier than computing generators of $\Bin(I)$.
A method for computing the monomial part is proposed in Tutorial~50 in~\cite{kreuzer2005computational}.

For the remainder of this section, we denote the set of indeterminates $\{x_1, \dots, x_n\}$ in~$P$ by~$X$, and let
$Y \subseteq X$ be a subset. Using the methods we developed so far, the binomial part of a cellular ideal can already
be partially determined.

\begin{remark}
    Let~$I$ be a $Y$-cellular ideal in~$P$.
    \begin{enumerate}[(a)]
        \item Since the indeterminates in $X\setminus Y$ are nilpotent modulo~$I$, only finitely many terms of
            $K[X\setminus Y]$ are not contained in~$I$. The fact that $I \cap K[Y]$ is saturated with respect to
            the product of all indeterminates in $Y$ implies that all monomials in $I$ have to be in
            $I \cap K[X \setminus Y]$. Hence the binomial part of $I \cap K[X \setminus Y]$ is generated by the
            monomial part of~$I$ and the binomials in the finitely many terms in $K[X \setminus Y]$ not contained
            in $\Mon(I)$.
        \item By assumption $I \cap K[Y]$ is saturated with respect to the product of indeterminates in $Y$. The
            binomial part of $I \cap K[Y]$ can therefore be computed using Corollary~\ref{cor:binomial_part_saturated}.
    \end{enumerate}
\end{remark}

An obstruction to computing the entire binomial part of~$I$ are binomials of the form $su- a vt \in I$, where $s,t$
are terms in $K[X \setminus Y]$, where $u,v$ are terms in $K[Y]$, and where $a\in K$.

\begin{definition}
Let $Y \subseteq X$, let $I$ be an ideal in~$P$ with $I : \langle \prod_{x_i \in Y} x_i \rangle = I$, and let $s,t$ be
fixed terms in $K[X \setminus Y]$. The ideal generated by all binomials of the form $s u - a vt$ with $a\in K$ and
terms $u,v \in K[Y]$ is called the \textbf{$(s,t)$-binomial part} of~$I$ and is denoted by $\Bin_{s,t}(I)$.
\end{definition}

If $s,t \notin I$ then the binomials in $\Bin_{s,t}(I)$ cannot be computed with any of the methods described above.

    \begin{example}\label{example:obstruction}
        Consider the ideal $I = \langle x^4,  y^4, x^2z^4 +xyz^2 +y^2,  x^3z^2 -x^3 -y^3 \rangle$ in
        $\mathbb{Q}[x,y,z]$. It is $Y$-cellular for $Y = \{z\}$. We have
        $\Mon(I) = \langle y^4,  xy^3,  x^2y^2,  x^3y,  x^4 \rangle$ and $I \cap \mathbb{Q}[z] = \langle 0 \rangle$.
        To compute $\Bin(I)$ we therefore need to search for binomials of the form $sz^a-\lambda z^bt$ with
        $a,b \in \mathbb{N}$, $\lambda \in K$ and $s,t \in \mathbb{T}_{x,y} \setminus \Mon(I)$, where $\mathbb{T}_{x,y}$
        is the monoid of all terms in $\mathbb{Q}[x,y]$ and $\Mon(I)$ is considered as a monoideal in $\mathbb{T}_{x,y}$.
    \end{example}

    For a set of indeterminates $Y \subseteq X$ we denote the localization with respect to the multiplicative set
    generated by $\prod_{x_i \in Y} x_i$ by $K[X]_Y$. When considering an ideal in $K[X]$ which is saturated with
    respect to $\prod_{x_i \in Y} x_i$ we can assume that the indeterminates in~$Y$ are invertible by
    Proposition~\ref{prop:invertible}. Consequently it is enough to look at ideals in~$K[X]_Y$.
    The next lemma shows that $\Bin_{s,t}(I)$ has an affine structure.

\begin{lemma}\label{lemma:binomial_structure}
        Let $Y \subseteq X$ with $\#Y = m$, and let $I \subseteq K[X]_Y$ be an ideal. Let $s$ and $t$ be terms in
        $ K[X \setminus Y]$.
\begin{enumerate}
\item[(a)] If $sY^a- \lambda t$ and $sY^b- \mu t$ are in $I$ for some $a, b \in \mathbb{Z}^m$ and
                $\lambda, \mu \in K^\times$, then we have $s Y^{a+kc} - (\lambda^{k+1}/\mu^k) t \in I$ for
                $c = a-b$ and all $k \in \mathbb{Z}$.

\item[(b)] Let $u_1, \dots u_s \in \mathbb{Z}^m$ and $\lambda_{u_1}, \dots, \lambda_{u_s} \in K^\times$ such that
                $sY^{u_1}-\lambda_{u_1}t \in I$ and $sY^{u_1+u_i}- \lambda_{u_1}\lambda_{u_i}t \in I$ for all
                $i=2, \dots, s$. Then we have $sY^w-\lambda_w t \in I$ for all
                $w = u_1 + k_2 u_2 + \cdots + k_s u_s$ with $k_2, \dots, k_s \in \mathbb{Z}$ and
                $\lambda_w = \lambda_{u_1}\lambda_{u_2}^{k_2} \cdots \lambda_{u_s}^{k_s}$.
\end{enumerate}
\end{lemma}
    
    \begin{proof}
        The following equation proves Part~(a).
        \begin{align*}
            sY^{a+k c}-(\lambda^{k+1}/\mu^k) t = (Y^{k c} + (\lambda/\mu)Y^{(k-1)c}+ \cdots + (\lambda/\mu)^c)(sY^a -\lambda t) \\
            - ((\lambda/\mu)Y^{k c} + (\lambda/\mu)^2Y^{(k-1)c} + \cdots + (\lambda/\mu)^k Y^c)(sY^b-\mu t)
        \end{align*}
        Part~(b) then follows by a direct calculation from Part~(a).
    \end{proof}

    Let $I \subseteq K[X]$ be an ideal which is saturated with respect to the product of the indeterminates in
    $Y \subseteq X$, and let $s,t$ be terms which are not contained in~$I$.
    The idea of the following algorithm is to pass to the ring $Q = K[X]_Y$ and search for elements of the form
    $(1,\lambda Y^a)$ with $a \in \mathbb{Z}^m$ and $\lambda \in K$ in the syzygy module~$\Syz(s,t)_{Q/QI}$.
    The $(s,t)$-binomial part is then determined by a unit lattice modulo the ideal $(IQ:\langle t \rangle) \cap K[Y]_Y$.
    This lattice and its associated character can be computed using Algorithm~\ref{alg:unit_lattice}. Note that
    this algorithm can also be applied to ideals in a localized polynomial ring, see Remark~\ref{remark:localization}.

\begin{algorithm}{\textbf{(Computing (s,t)-Binomial Parts)}}\label{alg:binomials_fixed}\\
Let $Y \subseteq X$, let $I \subseteq K[X]$ be an ideal with $I : \langle \prod_{x_i \in Y} x_i \rangle = I$,
and assume that $Y = \{x_1, \dots, x_m\}$. Let $s,t \in K[X \setminus Y]$ be terms with $s,t \notin I$. Consider
the following instructions.
\begin{enumerate}
\item[(1)] Form the ring $Q = K[X]_Y$.

\item[(2)] Compute generators $(f_1, g_1), \dots, (f_k, g_k)$ of
                $S \coloneqq \Syz(s,t)_{Q/IQ} \cap K[Y]_Y$.

\item[(3)] If $\langle f_1, \dots, f_k \rangle = \langle 1 \rangle$, then compute $h \in K[Y]_Y$ such that $(1, h) \in S$.
                If $\langle f_1, \dots, f_k \rangle \neq \langle 1 \rangle$ or $h$ is not a unit
                in $Q/(IQ: \langle t \rangle)$, return the zero ideal.

\item[(4)] Using Algorithm~\ref{alg:unit_lattice} compute the unit lattice
                $\Lambda= \langle v_1, \dots, v_r \rangle \subseteq \mathbb{Z}^{m+1}$ and the associated
                character $\rho$ of $(x_1, \dots, x_m, -h)$ modulo $(IQ:\langle t \rangle) \cap K[Y]_Y$.

\item[(5)] Consider the equation over $\mathbb{Z}$ in the indeterminates $z_1, \dots, z_r$ given by
                $$z_1 v_{1,m+1} + \cdots + z_r v_{r, m+1} = 1.$$
                Compute $u_1, \dots, u_s \in \mathbb{Z}^r$ such that all integer solutions of this equation
                are given by $u_1+\mathbb{Z}u_2 + \cdots + \mathbb{Z}u_s$.

\item[(6)] For $i=1, \dots, r$ let $v'_i$ be the tuple consisting of the first $m$ components of $v_i$. Let
                $M \in \Mat_{m, r}(\mathbb{Z})$ be the matrix whose columns are given by $v'_1, \dots, v'_r$, and form
                the ideal
                $$ J = \left\langle sY^{Mv} - \rho(Mv) t \mid v \in \{u_1, u_1+u_2, \dots, u_1+u_s\} \right\rangle
                \subseteq Q.$$

\item[(7)] Return $J \cap K[X]$.
\end{enumerate}
This is an algorithm which computes $\Bin_{s,t}(I)$.
\end{algorithm}

\begin{proof}
        Assume that $sY^{a^+}-\lambda Y^{a^-}t \in \Bin_{s,t}(I)$ for $\lambda \in K^\times$ and
        $a = (a_1, \dots, a_m)$ in $\mathbb{Z}^m$. Then we have $s - \lambda Y^{-a} t \in IQ$ which shows
        $(1, -\lambda Y^{-a}) \in S$. In Step~(3) we therefore have $\langle f_1, \dots, f_k \rangle = \langle 1 \rangle$.
        Hence, there exists $h \in K[Y]_Y$ such that $(1, h) \in S$. Since $(1, h) \in S \subseteq \Syz(s,t)$, we have
        $s + ht \in IQ$ and therefore $$ -Y^a(s + ht)+(sY^a - \lambda t) = - h Y^a t - \lambda t \in IQ.$$
        This shows $- h Y^a -\lambda \in (IQ:\langle t \rangle)$, proving that $h$ is a unit in
        $Q/(IQ : \langle t \rangle)$. It further shows that we have $v = (a_1, \dots, a_m ,1) \in \Lambda$ with
        $\rho(v) = \lambda$. Therefore there exists $w \in u_1+\mathbb{Z}u_2 + \cdots + \mathbb{Z}u_s$ such that
        $Mw = (a_1, \dots, a_m)$. Now Lemma~\ref{lemma:binomial_structure} shows $sY^a-\rho(Mw) t\in J$. Finally,
        $t \notin I$ implies $\rho(Mw) = \lambda$, and we get $sY^{a^+}-\rho(Mw)Y^{a^-} t \in J \cap K[X]$.
\end{proof}

    \begin{example}
        Let $I \subseteq \mathbb{Q}[x,y,z]$ be given as in Example~\ref{example:obstruction}. We choose $s=x^3$ and
        $t=y^3$ from $\mathbb{T}_{x,y} \setminus \Mon(I)$. The module $S$ in Step~(1) is generated by
        $(z^2-1,-1)$ and $(1, z^4+z^2)$. This yields $h = z^4+z^2$. Next, we compute
        $(I: \langle t \rangle) = \langle z^6-z^2+1 \rangle$, and confirm that $h$ is a unit
        in $\mathbb{Q}[z]/\langle z^6-z^2+1 \rangle$. The unit lattice of $(z, -h)$ modulo
        $\langle z^6-z^2+1 \rangle$ is given by $\Lambda = \langle (6, 1), (0,16) \rangle$ and the associated
        character is $\rho : \Lambda \rightarrow \mathbb{Q}^\times$ defined by $\rho((6, 1)) = \rho((0,16)) = 1$.
        From the first generator of $\Lambda$ we immediately obtain $x^3z^6-y^3 \in I$.
    \end{example}

    The preceding algorithm allows us to compute the binomial part of a single cellular ideal. But to obtain the
    binomial part of an intersection of cellular ideals, we can not simply intersect their binomial parts. This is
    because the intersection of binomial ideals is in general not binomial.

    Instead we use the following observation.
    Let $I = I_1 \cap \cdots \cap I_k$ be a decomposition into $Y_i$ cellular ideals. Then for every binomial $f$ in $I$
    there is a intersection $Y$ of elements from $\{Y_1, \dots, Y_k\}$ such that
    $f \in \bigcap_{Y \not\subseteq Y_i} \Mon(I_i)$ and $f \in \Bin_{s,t}(\bigcap_{Y \subseteq Y_i} I_i)$ for some
    terms $s,t$ in $K[X \setminus Y]$.

\begin{algorithm}{\textbf{(Computing Binomial Parts)}}\label{alg:binomial_part}\\
Let $I \subseteq K[X]$ be an ideal. Consider the following instructions.
\begin{enumerate}[(1)]
            \item Let $B = [\;]$.
            \item Using Algorithm~\ref{alg:cellular_decomposition} compute a decomposition $I = I_1 \cap \cdots \cap I_k$
                where $I_i$ is a $Y_i$-cellular ideal for some $Y_i \subseteq X$.
            \item For all $x_i \in X\setminus (\bigcap_{i=1}^k Y_i)$ let $\delta_i$ be minimal
                such that $x_i^{\delta_i} \in I_j$ for all~$j$ with $x_i \in X \setminus Y_j$.
            \item For all sets $S \subseteq \{Y_1, \dots, Y_k\}$ let $Y = \bigcap_{s \in S} s$
                and perform steps (5)--(8). Note, that this includes the empty intersection given by $X$.
            \item Compute $M_Y = \bigcap_{Y \not\subseteq Y_i} \Mon(I_i)$ and
                $J_Y = \bigcap_{Y \subseteq Y_i} I_i$.
            \item Compute the finite set $T_Y$ of terms in $K[X \setminus Y]$ such that the exponent of $x_i$ is
                smaller than $\delta_i$ for all $x_i \in X \setminus Y$.
            \item For each pair $s,t \in T_Y$ compute $\Bin_{s,t}(J_Y)$. If $s \in J_Y$ or $t \in J_Y$, then it is given
                by $\langle s \rangle$, $\langle t \rangle$ or $\langle s, t \rangle$. Otherwise we can use
                Algorithm~\ref{alg:binomials_fixed}.
            \item Add the generators of $M_Y \cap \sum_{s,t \in T_Y} \Bin_{s,t}(J_Y)$ to $B$.
            \item Return $B$.
\end{enumerate}
        This is an algorithm which computes generators of the binomial part $\Bin(I)$.
\end{algorithm}

\begin{proof}
        The set $B$ consists of binomials since the intersection of a monomial ideal and a binomial ideal in Step~(8)
        is again binomial. The containment $\langle B \rangle \subseteq \Bin(I)$ is clear, since
        $M_Y \cap J_Y \subseteq I$ for all $Y \subseteq X$.

        To show the opposite inclusion assume that
        $f = x_1^{a_1} \cdots x_n^{a_n} - \lambda x_1^{b_1} \cdots x_n^{b_n} \in I$ where $a_i, b_i \in \mathbb{N}$ and
        $\lambda \in K$. Let
        $$Y = \{x_i \mid a_i \geq \delta_i \text{ or } b_i \geq \delta_i \} \cup \bigcap_{i=1}^k Y_i.$$
        We now show that $f$ is contained in $M_Y \cap \sum_{s,t \in T_Y} \Bin_{s,t}(J_Y)$.
        For all $i=1, \dots, k$ with $Y \not\subseteq Y_i$ there exists $x_{j_i} \in Y$ which is not contained in $Y_i$.
        The indeterminate is therefore nilpotent modulo $I_i$ and we have $x_{j_i}^{\delta_{j_i}} \in I_i$. By definition
        of $Y$ we then have $x_{j_i}^{a_i} \in I_i$ or $x_{j_i}^{b_i} \in I_i$. This proves $f \in M_Y$.

        Since for all $I_i$ with $Y \subseteq Y_i$ we have $I_i : \langle \prod_{x_i \in Y} x_i \rangle = I_i$ the ideal
        $J_Y$ also satisfies $J_Y : \langle \prod_{x_i \in Y} x_i \rangle = J_Y$. For
        $s = \prod_{x_i \in X \setminus Y} x_i^{a_i}$ and $t = \prod_{x_i \in X \setminus Y} x_i^{b_i}$ we then have
        $s,t \in T_Y$ since $a_i < \delta_i$ and $b_i < \delta_i$. Therefore $f \in \Bin_{s,t}(J_Y)$ and hence $f$ is
        contained in the ideal $M_Y \cap \sum_{s,t \in T_Y} \Bin_{s,t}(J_Y)$.

        Now we need to show that $Y$ is of the form required in Step~(4). But if we choose
        $Y' = \bigcap_{Y \subseteq Y_i}$ then $J_Y = J_{Y'}$ and $M_Y = M_{Y'}$. It
        is therefore enough to only consider non-empty intersections of elements from $\{Y_1, \dots, Y_k\}$.
    \end{proof}

    The following remark provides some details on how to perform the steps of this algorithm.

    \begin{remark}\label{remark:bin_alg}
        \begin{enumerate}[(a)]
            \item In Step~(7), if $s$ and $t$ have a common factor $w \in \mathbb{T}^n$, it is enough to form
                the terms $s' = s/w$ and $t' = t/w$, compute $\Bin_{s', t'}(J_Y : \langle w \rangle)$, and then
                multiply its binomial generators with $w$ to obtain the ideal $\Bin_{s,t}(J_Y)$.
            \item If in Step~(6) we have $T_Y = \{1\}$, then $\Bin_{s,t}(J_Y) = \Bin(J \cap K[Y])$. This happens for
                example in the case $Y=X$. Since $J \cap K[Y]$ is saturated with respect to the product of all
                indeterminates in $Y$, we can compute $\Bin(J \cap K[Y])$ using Corollary~\ref{cor:binomial_part_saturated}.
            \item For $Y = \emptyset$, the binomials in Step~(7) can be computed by checking for every pair of
                terms $s,t \in T_Y$ whether $s-\lambda t \in I$ for some $\lambda \in K$. This can be achieved
                by checking whether the normal forms of $s$ and $t$ are scalar multiples of each other. In the final
                section of this paper we show some optimizations for this case.
        \end{enumerate}
    \end{remark}

    Let us conclude this section by applying the algorithm for computing the binomial part to a concrete example.

    \begin{example}
       Let $I = \langle x^3z^4 + x^2yz^2 + xy^2, x^2y + xy^2 + y^3 \rangle \subseteq \mathbb{Q}[x,y,z]$.
       We compute a cellular decomposition $I = I_1 \cap I_2 \cap I_3$ into the $Y_i$-cellular ideals
       $I_i$ given by
       \begin{align*}
           I_1 &= \langle z^6 -1,  yz^4 -xz^2 -yz^2 +x,  x^2 +xy +y^2 \rangle, &&Y_1 = \{x,y,z\},\\
           I_2 &= \langle z^4, y \rangle, &&Y_2 = \{x\}, \\
           I_3 &= \langle x^3z^4 + x^2yz^2 + xy^2, x^2y + xy^2 + y^3, x^5 \rangle, &&Y_3 = \{z\}.
       \end{align*}
       From this we obtain the bounds $\delta_x = 5$, $\delta_y = 6$ and $\delta_z = 4$. For each
       $Y$ in the set $\{\emptyset, \{x\}, \{x,y,z\}, \{z\}\}$ we then compute the ideals $J_Y$ and $M_Y$, the set
       $T_Y$ and the ideal $\mathcal{B}_Y = \sum_{s,t \in T_Y} \Bin_{s,t}(J_Y)$.
       For $Y = \emptyset$ we have
       \begin{align*}
           J_Y &= I_1 \cap I_2 \cap I_3, \\
           M_Y &= \langle 1 \rangle, \\
           T_Y &= \mathbb{T}_{x,y,z}\setminus \langle x^5, y^6, z^4 \rangle, \\
           \mathcal{B}_Y &= \langle x^3y - y^4 \rangle.
       \end{align*}
       For $Y = \{x\}$ we have
       \begin{align*}
           J_Y &=I_1 \cap I_2, \\
           M_Y &= \Mon(I_3) = \langle x^2y^3, x^5, y^6, xy^5, x^4y^2 \rangle, \\
           T_Y &= \mathbb{T}_{y,z} \setminus \langle y^6, z^4 \rangle, \\
           \mathcal{B}_Y &= \langle x^3y - y^4 \rangle.
       \end{align*}
       For $Y = \{x,y,z\}$ we have
       \begin{align*}
           J_Y &= I_1, \\
           M_Y &= \Mon(I_2) \cap \Mon(I_3) = \langle x^2y^3, y^6, xy^5, x^4y^2, x^5y, x^5z^4 \rangle.
       \end{align*}
       Since $J_Y$ is saturated with respect to the product of all indeterminates, we directly compute
       $\mathcal{B}_Y = \Bin(J_Y) = \langle x^3 - y^3, z^6 - 1 \rangle$.
       For $Y = \{z\}$ we have
       \begin{align*}
           J_Y &= I_1 \cap I_3, \\
           M_Y &= \Mon(I_2) = \langle y, z^4 \rangle, \\
           T_Y &= \mathbb{T}_{x,y} \setminus \langle x^5, y^6 \rangle, \\
           \mathcal{B}_Y &= \langle x^2y^3z^6 -x^2y^3,  y^5z^6 -y^5,  x^4z^6 -xy^3,  x^3y -y^4 \rangle.
       \end{align*}
       Altogether we obtain
       $$\Bin(I) = \langle y^5z^6 -y^5,  x^4z^6 -xy^3,  x^5z^4 -x^2y^3z^4,  x^3y -y^4 \rangle.$$
    \end{example}

\bigskip\bigbreak
%%%%%%%%%%%%%%%%%%%%%%%%%%%%%%%%%%%%%%%%%%%%%%%%%%%
%
%  Section 7: Optimizations
%
%%%%%%%%%%%%%%%%%%%%%%%%%%%%%%%%%%%%%%%%%%%%%%%%%%%

\section{Optimizations}

In the final section of this paper we show two ways in which Algorithm~\ref{alg:binomial_part} can be optimized.

%%%%%%%%%%%%%%%%%%%%%%%%%%%%%%%%%%%%%
% Subsection 7.1
%%%%%%%%%%%%%%%%%%%%%%%%%%%%%%%%%%%%%
\medskip
\subsection*{Binomials in Vector Spaces}

If $Y = \emptyset$, then $\Bin_{s,t}(J_Y)$ in Step~(7) of Algorithm~\ref{alg:binomial_part}
can be computed by checking whether there exists $\lambda \in K$ such that $s-\lambda t \in I$. This needs to be
done for all terms $s,t \in T_Y$. In other words we need to determine the ideal generated by all binomials in
$\Bin(I)$ whose support is contained in $T_Y$. Instead of checking every pair of terms individually, we can use the
following method. To describe it, we translate the concept of a binomial to vector spaces.

\begin{definition}
Let $V$ be a finite-dimensional vector space over a field $K$ with basis $e_1, \dots, e_n$, 
and let $U \subseteq V$ be a subspace.
\begin{enumerate}
\item[(a)] A \textbf{binomial} in $U$ is an element of the form $c_i e_i - c_j e_j$ with
                $c_i, c_j \in K$ and $i \neq j$.

\item[(b)] The subspace spanned by all binomials contained in $U$ is denoted by $\Bin(U)$. We call it the
                \textbf{binomial part} of $U$.
\end{enumerate}
\end{definition}

Without loss of generality, we can confine our search to binomials of the form $e_i - c e_j$.

\begin{algorithm}{\textbf{(Computing the Binomial Part of a Vector Subspace)}}\label{alg:vr_binomials}\\
        Let $V$ be a finite-dimensional vector space over a field $K$ with basis $e_1, \dots, e_n$,
        let $U \subseteq V$ be a subspace spanned by $b_1, \dots, b_k$, and let the coordinates of $b_i$ with respect
        to $(e_1, \dots, e_n)$ be given by $c_i \in K^n$. The following steps define an algorithm which computes
        generators of the binomial part $\Bin(U)$.
\begin{enumerate}[(1)]
\item[(1)] Let $S = [\;]$.

\item[(2)] Compute the reduced row echelon form $M$ of the matrix
$$
(c_1, \dots, c_k)^{tr} \in \Mat_{k \times n}(K).
$$

\item[(3)] If a row of $M$ has at most two non-zero components, add the corresponding binomial to~$S$.

\item[(4)] For each pair of rows $(m_i, m_j)$ of $M$ check if there exists $c\in K^\times$ such that $r_i-c r_j$
                has at most two non-zero components. Add the corresponding binomials to $S$.

\item[(5)] Return $S$.
\end{enumerate}
\end{algorithm}

\begin{proof}
        The coordinates of an element in $U$ are given by a linear combination of the rows of $M$.
        Since $M$ is in reduced row echelon form, the sum of more than two rows already corresponds to an element
        of~$U$ with at least three non-zero coordinates. Hence we conclude that the coordinates of a binomial can only
        be given by a linear combination of at most two rows.
\end{proof}

An alternative method for computing the binomial part of a vector subspace using matroid theory is described
in~\cite{kahle2022short}. We can now compute the binomial part restricted to a finite set of terms as follows.

\begin{algorithm}{\textbf{(The Binomial Part Restricted to a Finite Set of Terms)}}\\\label{alg:restricted_binomial_part}
Let $I$ be an ideal, $\sigma$ a term ordering, and $T \subseteq \mathbb{T}^n$ a finite set of terms.
Assume that $T$ is closed with respect to $\sigma$, i.e., if $t \in T$ then for all $s \in \mathbb{T}^n$ with
$s \leq_\sigma t$ we have $s \in T$. The following instructions form an algorithm which computes the ideal
generated by all binomials $f \in I$ with $\Supp(f) \subseteq T$.
\begin{enumerate}
\item[(1)] Let $S = [\;]$.

\item[(2)] Compute a $\sigma$-Gröbner basis $G$ of $I$.

\item[(3)] Add all $g \in G$ with $\Supp(g) \subseteq T$ to $S$.

\item[(4)] For every term $t \in T$ if $t \notin \LT_\sigma(S)$ and there exists $h \in S$ and
                $s \in \mathbb{T}^n$ such that $\LT_\sigma(sh) = t$, then add $sh$ to $S$.

\item[(5)] Let $V$ be the vector space over~$K$ spanned by~$T$, and let $U = \vspan_K(S)$ be
the subspace of~$V$ spanned by~$U$. Apply Algorithm~\ref{alg:vr_binomials} to compute $\Bin(U)$ and return it.

\end{enumerate}
\end{algorithm}

\begin{proof}
        Clearly, the algorithm terminates since $T$ is finite. Let us now show that after performing Steps
        (1) to (4) the set $S$ generates the subspace of $V$ spanned by all polynomials $f \in I$ with
        $\Supp(f) \subseteq T$. We denote this subspace by $U$. It is easy to see that $\vspan_K(S)$ is
        contained in $U$. Suppose there exists $f \in U$ such that $f \notin \vspan_K(S)$. Then we can choose $f$
        such that it has a $\sigma$-minimal leading term. Since $f \in I$ and $T$ is closed with respect to $\sigma$,
        there exist $g \in G$ with $\Supp(g) \subseteq T$ and $t \in \mathbb{T}^n$ such that
        $\LT_\sigma(f) = \LT_\sigma(gt)$. It follows from the condition in Step~(4) that there exists $h \in S$ with
        $\LT_\sigma(f) = \LT_\sigma(h)$. Since $f \notin \vspan_K(S)$, we have $f-h \notin \vspan_K(S)$. But
        $\LT_\sigma(f-h)$ is smaller than $\LT_\sigma(f)$. This is a contradiction.
\end{proof}

This algorithm can be used as an optimization in Algorithm~\ref{alg:binomial_part}.

\begin{corollary}
In Algorithm~\ref{alg:binomial_part} if in Step~(4) we choose $Y = \emptyset$, then replace Steps~(5) to~(8) by
the following steps.
\begin{enumerate}
\item[(5')] We have $J_Y = I$ and $M_Y = \langle 1 \rangle$. This step can therefore be omitted.

\item[(6')] Compute the set $T \subseteq \mathbb{T}^n$ of terms such that the exponent of $x_i$ is
                smaller than $\delta_i$ for $i=1, \dots, n$.

\item[(7')] Using Algorithm~\ref{alg:restricted_binomial_part} compute the ideal $\Bin_T(I)$ generated by
                all binomials $f \in I$ with $\Supp(f) \subseteq T$.

\item[(8')] Add the generators of $\Bin_T(I)$ to $B$.
\end{enumerate}
The result is an algorithm which computes the binomial part of~$I$.
\end{corollary}

Algorithm~\ref{alg:restricted_binomial_part} can also be used to compute the binomial part up to a degree bound.

\begin{remark}
Let $\delta \geq 1$. If we choose a degree compatible term ordering and let $T \subseteq \mathbb{T}^n$ be
the set of all terms~$t$ with $\deg(t) \le \delta$, then Algorithm~\ref{alg:restricted_binomial_part} can be
used to compute the binomial part up to the degree bound~$\delta$.

However, note that in general this can not be used to compute all of~$\Bin(I)$, since no degree bound for the
generators of~$\Bin(I)$ is known.
\end{remark}

%%%%%%%%%%%%%%%%%%%%%%%%%%%%%%%%%%%%%%
% Subsection 7.2
%%%%%%%%%%%%%%%%%%%%%%%%%%%%%%%%%%%%%%
\bigbreak
\subsection*{Binomial Parts of Radical Ideals}
   
In the remainder of this section we show that Algorithm~\ref{alg:binomial_part} can be simplified 
if the ideal $I$ is radical.

\begin{definition}
Let $Y \subseteq X$. A binomial ideal $I$ in $K[X]$ is called $Y$-\textbf{mesoprime} 
if the following conditions hold.
\begin{enumerate}
\item[(a)] $I = I : \langle \prod_{x_i \in Y} x_i \rangle^\infty$.

\item[(b)] $x_i \in I$ for all $x_i \in X \setminus Y$.
\end{enumerate}
\end{definition}

Note that every $Y$-mesoprime ideal, is $Y$-cellular. Also every $Y$-mesoprime ideal is of the form
$J + \langle X \setminus Y \rangle$ where $J$ is a lattice ideal in $K[Y]$.

\begin{lemma}\label{lemma:binomial_prime}
        The binomial part of a prime ideal is mesoprime.
\end{lemma}

\begin{proof}
        Let $\mathfrak{p}$ be a prime ideal. If $\mathfrak{p} : \langle x_i \rangle \neq \mathfrak{p}$, then there
        exists $f \in K[X] \setminus \mathfrak{p}$ such that $x_i f \in \mathfrak{p}$, which implies $x_i \in \mathfrak{p}$.
\end{proof}

The following lemma implies that for an ideal whose binomial part is mesoprime it is not necessary to consider
$(s,t)$-binomial parts.

\begin{lemma}\label{lemma:bin_cellular}
        Let $I \subseteq K[X]$ be a proper, $Y$-cellular ideal for some $Y \subseteq X$. Then $I$ does not contain
        binomials of the form $t-\lambda s$ with $\lambda \in K$ and terms $t \in K[Y]$ and $s \notin K[Y]$.
\end{lemma}

\begin{proof}
        To prove (a), suppose that a binomial $f = t-\lambda s$ as above is contained in~$I$. Then $x_i$ divides $s$
        for some $x_i \notin Y$. Since $I$ is $Y$-cellular we have $x_i^d \in I$ for some $d \in \mathbb{N}$.
        The equation
        $$ t^3+\lambda^3 s^3 = t^2 f + \lambda^2s^2 f + \lambda ts f $$ shows
        that $t^\ell + \lambda^\ell s^\ell \in I$ for some $\ell \geq d$. This implies $t^\ell \in I \cap K[Y]$,
        which contradicts the assumption that $I$ is $Y$-cellular.
\end{proof}

If the binomial parts of ideals are mesoprime, then the binomial part of their intersection can be obtained as follows.

\begin{corollary}\label{cor:mesoprime_intersection}
        Let $I_1, \dots, I_k \subseteq K[X]$ be ideals whose binomial parts $\Bin(I_i)$ are $Y_i$-mesoprime for some
        $Y_i \subseteq X$. Assume that $\Bin(I_i) = J_i + \langle X \setminus Y_i \rangle$ where $J_i$
        is a lattice ideal in $K[Y]$. Then we have
        \[\Bin \left(I_1 \cap \cdots \cap I_k \right) = \sum_{S \subseteq \{1, \dots, k\}} \bigcap_{j \notin S}
        \langle X \setminus Y_j \rangle \cap J_S \]
        where $J_S$ is the lattice ideal $\Bin(\bigcap_{i \in S} J_i)$. For $S = \emptyset$, we use $I_S = \langle 1 \rangle$ here.
\end{corollary}

\begin{proof}
        The right hand side is a binomial ideal since the intersection of a monomial ideal and a binomial ideal is again
        binomial. Also, it is easy to see that the right hand side is contained in the left hand side. To show the
        opposite inclusion let $f$ be a binomial in $I_1 \cap \cdots \cap I_k$. By Lemma~\ref{lemma:bin_cellular} we
        have for each $i=1, \dots, k$ either $f \in K[Y_i]$ or $f \in \langle X \setminus Y_i \rangle$.
        Thus, there exists $S \subseteq \{1, \dots, k\}$ such that $f \in K[Y_i]$ for all $i \in S$ and
        $f \in \langle X \setminus Y_i \rangle$ for all $i \notin S$. This means $f$ is contained
        in $\bigcap_{j \notin S} \cap J_S$.
\end{proof}

The binomial parts of the form $\Bin(\bigcap_{i \in S} J_i)$ in the corollary can be computed using
Algorithm~\ref{alg:lattice_intersection} and Proposition~\ref{prop:lattice_to_binomial}. 
Now obtain the following algorithm for computing the binomial part of a radical ideal.

\begin{algorithm}{\bf (Computing the Binomial Part of a Radical Ideal)}\label{alg:bin_of_rad}\\
Let $I \subseteq P$ be a radical ideal. Consider the following sequence of instructions.
\begin{enumerate}
\item[(1)] Compute the prime decomposition $I = \mathfrak{p}_1 \cap \cdots \cap \mathfrak{p}_k$.

\item[(2)] For each $\mathfrak{p}_i$ compute the maximal set $Y_i \subseteq X$ such that we have
                $\mathfrak{p}_i : \langle \prod_{x_j \in Y} x_j \rangle = \mathfrak{p}_i$.

\item[(3)] For $i = 1, \dots, k$ compute $\Bin(\mathfrak{p}_i \cap K[Y_i]) = J_i$, and obtain
                $\Bin(\mathfrak{p}_i) = J_i + \langle X \setminus Y_i \rangle$.

\item[(4)] Return the ideal
$\sum_{S \subseteq \{1, \dots, k\}} \bigcap_{j \notin S} \langle X \setminus Y_j \rangle \cap J_S$.
\end{enumerate}
This is an algorithm which computes the binomial part $\Bin(I)$.
\end{algorithm}

\begin{proof}
The binomial part of $\mathfrak{p}_i$ is mesoprime by Lemma~\ref{lemma:binomial_prime}.
Now Corollary~\ref{cor:mesoprime_intersection} shows that Step~(4) returns $\Bin(I)$.
\end{proof}

\bigbreak

\bibliographystyle{acm}
\bibliography{binomial_part}

\end{document}